\documentclass[12pt]{article}

\usepackage{subcaption}

\usepackage{float}

\usepackage{multirow}

\usepackage{tikz}

\usetikzlibrary{calc,backgrounds,arrows,matrix}

\usepackage{enumerate}

\usepackage{color}
\definecolor{lightblue}{rgb}{0,0.2,0.5}

\usepackage[colorlinks=true, urlcolor=blue,linkcolor=blue, citecolor=lightblue]{hyperref}

\usepackage[round,sort,comma,numbers]{natbib}

\bibpunct{\textcolor{lightblue}{(}}{\textcolor{lightblue}{)}}{,}{a}{}{;}

\usepackage{amssymb,amsmath}

\usepackage{graphicx}

\usepackage[ruled,vlined]{algorithm2e}

\usepackage{color}
\usepackage{mathrsfs}

\usepackage{mathtools}

\usepackage{amsthm}
\usepackage{nameref}
\usepackage{bigints}

\DeclareMathAlphabet{\eufrak}{U}{}{}{}
\DeclarePairedDelimiter\floor{\lfloor}{\rfloor}

\SetMathAlphabet\eufrak{normal}{U}{euf}{m}{n}
\SetMathAlphabet\eufrak{bold}{U}{euf}{b}{n}

\allowdisplaybreaks

 \def\qu{{\mathord{\mathbb Z}}}

 \def\sZZ{{\rm Z\kern-.45em{}Z}}

 \def\sQQ{{\kern 0.27em \vrule height1.45ex width0.03em depth0em
           \kern-0.30em \rm Q}}
 \def\qu{{\mathchoice
         {\sQQ}
         {\sQQ}
   {\kern 0.225em \vrule height1.05ex width0.025em depth0em \kern-0.25em \rm Q}
   {\kern 0.180em \vrule height0.78ex width0.020em depth0em \kern-0.20em \rm Q}
         }}
 \def\sGG{{\kern 0.27em \vrule height1.45ex width0.03em depth0em
           \kern-0.30em \rm G}}
 \def\gg{{\mathchoice
         {\sGG}
         {\sGG}
   {\kern 0.225em \vrule height1.05ex width0.025em depth0em \kern-0.25em \rm G}
   {\kern 0.180em \vrule height0.78ex width0.020em depth0em \kern-0.20em \rm G}
         }}

 \newtheorem{prop}{Proposition}[section]
 \newtheorem{lemma}[prop]{Lemma}
 \newtheorem{definition}[prop]{Definition}
 \newtheorem{corollary}[prop]{Corollary}
 \newtheorem{theorem}[prop]{Theorem}
 \newtheorem{remark}[prop]{Remark}
 
 \newtheorem{assumption}{Assumption}

 \newtheoremstyle{named}{}{}{\itshape}{}{\bfseries}{.}{.5em}{#1 \thmnote{#3}}
\theoremstyle{named}

\numberwithin{equation}{section}

\newcommand{\abs}[1]{\lvert#1\rvert}
\newcommand{\norm}[1]{\left\lVert#1\right\rVert}

 \def\P{{\mathord{\mathbb P}}}

\newcommand{\re}{\mathrm{e}}

 \newcounter{hyp}
\usepackage{aeguill}

\newenvironment{Proof}{\removelastskip\par\medskip \noindent{\em Proof.} \rm}{\penalty-20\null\hfill$\square$\par\medbreak}

\def\bprf{\begin{Proof}}
\def\nprf{\end{Proof}}
\def\bdes{\begin{description}}
\def\ndes{\end{description}}

\newtheorem{thm}{Theorem}[section]

\def\bdef{\begin{defn}}
\def\ndef{\end{defn}}
\def\bthm{\begin{thm}}
\def\nthm{\end{thm}}
\def\bprop{\begin{prop}}
\def\nprop{\end{prop}}
\def\brmk{\begin{remark}}
\def\nrmk{\end{remark}}
\def\bexa{\begin{exa}}
\def\nexa{\end{exa}}
\def\blem{\begin{lem}}
\def\nlem{\end{lem}}
\def\bcor{\begin{cor}}
\def\ncor{\end{cor}}
\def\bexe{\begin{exe}}
\def\nexe{\end{exe}}

\newcommand{\E}{\mathbb{E}}

\newcommand{\real}{\mathbb{R}}

\def\og{\leavevmode\raise.3ex
     \hbox{$\scriptscriptstyle\langle\!\langle$~}}
\def\fg{\leavevmode\raise.3ex
     \hbox{~$\!\scriptscriptstyle\,\rangle\!\rangle$}~}

\title{\Huge
A constructive approach to existence of equilibria in time-inconsistent stochastic control problems
}

\author{
 Jiang Yu Nguwi\footnote{\href{mailto:nguw0003@e.ntu.edu.sg}{nguw0003@e.ntu.edu.sg}
 }
 \qquad Nicolas Privault\footnote{
\href{mailto:nprivault@ntu.edu.sg}{nprivault@ntu.edu.sg}
 }
 \\
\small
Division of Mathematical Sciences
\\
\small
School of Physical and Mathematical Sciences
\\
\small
Nanyang Technological University
\\
\small
21 Nanyang Link, Singapore 637371
}

\allowdisplaybreaks

\usepackage{empheq}

\usepackage{bbm}

\makeatletter
\newcommand*\rel@kern[1]{\kern#1\dimexpr\macc@kerna}
\newcommand*\widebar[1]{
  \begingroup
  \def\mathaccent##1##2{
    \rel@kern{0.8}
    \overline{\rel@kern{-0.8}\macc@nucleus\rel@kern{0.2}}
    \rel@kern{-0.2}
  }
  \macc@depth\@ne
  \let\math@bgroup\@empty \let\math@egroup\macc@set@skewchar
  \mathsurround\z@ \frozen@everymath{\mathgroup\macc@group\relax}
  \macc@set@skewchar\relax
  \let\mathaccentV\macc@nested@a
  \macc@nested@a\relax111{#1}
  \endgroup
}
\makeatother

\let\oldcitet=\citet
\let\oldcitep=\citep
\renewcommand{\cite}[1]{\textcolor[rgb]{0,0,1}{\oldcitet{#1}}}
\renewcommand{\citet}[1]{\textcolor[rgb]{0,0,1}{\oldcitet{#1}}}
\renewcommand{\citep}[1]{\textcolor[rgb]{0,0,1}{\oldcitep{#1}}}

\begin{document}
\maketitle

\baselineskip0.6cm

\vspace{-0.6cm}

\begin{abstract}
    We extend the construction of equilibria
    for linear-quadratic and
    mean-variance portfolio problems
    available in the literature
    to a large class of mean-field time-inconsistent
    stochastic control problems
    in continuous time.
Our approach relies on a time discretization of the control problem via $n$-person games, which are characterized via the maximum principle using Backward Stochastic Differential Equations (BSDEs). The existence of equilibria is proved by applying weak convergence arguments to the solutions of $n$-person games. A numerical implementation is provided by approximating $n$-person games using finite Markov chains.
\end{abstract}



\baselineskip0.7cm

\parskip-0.1cm

\section{Introduction}
Stochastic control theory aims at
 optimizing a
    time-dependent
functional
    parameterized
by a controlled random state process,
with
    applications
    to numerous problems in physics, 
     biology, finance, economics, etc.
For this, the most commonly used approaches rely on
 Pontryagin's maximum principle and on
 Hamilton-Jacobi-Bellman (HJB) equations,
 see e.g. \cite{yong1999stochastic} and \cite{fleming2006controlled}
 for classical results on stochastic control theory.
     This approach deals with time-consistent stochastic control problems,
     in which an optimal strategy today remains optimal in the future.

\medskip

    However, many stochastic control problems are time-inconsistent
    in the sense that an optimal strategy today may not be optimal in the future.
    This is the case for example in the framework of a
    production economy with time-varying preferences,
    or in the nonlinear setting of mean-variance
    portfolio optimization, which cannot be directly treated using
    the dynamic programming principle and HJB equations.
    Such problems have recently been the
    object of increased attention, see e.g.\
    \cite{bjork2014theory} and \cite{bjork2017time}.

\medskip

There are two common formulations for time-inconsistent problems.
The first approach is to fix an initial time, to solve the problem given this
initial time, and to stick to this pre-committed optimal policy for the
remaining time.
See for example \cite{zhou2000continuous} for the solution
of mean-variance portfolio selection problem using pre-committed strategies.

\medskip

The second approach,
introduced by \cite{ekeland2006being} in the deterministic setting,
is to formulate time-inconsistent problems in a game-theoretic setting
using equilibrium controls.
This approach, which uses an HJB-type equation
to characterize the equilibrium controls,
has been extended in
 \cite{bjork2014theory} and \cite{bjork2017time}
 to stochastic mean-field control problems
  in both discrete and continuous time.
  In \cite{hu2012time,hu2017time},
a related characterization has been proposed
by the maximum principle in a linear-quadratic model,
where the SDE is linear and the mean-field objective functional is quadratic.
This characterization argument has been later
extended to general mean-field objective functionals
in \cite{djehiche2016characterization}.

\medskip

However, no general results are available on the
existence of equilibrium controls,
except
in special cases such as the linear-quadratic model of \cite{hu2012time}.
In addition,
no numerical construction of equilibrium controls has been provided so far,
 except in mean-variance portfolio selection, see \cite{wang2011continuous}.

 \medskip

 In this paper, we present a constructive approach to the existence of
 equilibrium controls for
a class of
mean-field time-inconsistent control problems,
together with its numerical implementation.
Let $(W_t)_{t\in [0,T]}$ denote a standard Brownian motion
generating the filtration $({\cal F}_t)_{t\in [0,T]}$.
Our results apply to the class of cost functionals of the form
\begin{equation}
\label{cost_functional_relaxed}
 J(t,\xi,\mu) = \mathbb{E}_t\left[
   g\big(X^{\xi,\mu}_{t,T}, \mathbb{E}_t \big[\Psi\big( X^{\xi,\mu}_{t,T}\big)\big]
   \big)
   +
   \int_t^T \int_U h\big(s, X^{\xi,\mu}_{t,s}, \mathbb{E}_t\big[\Phi\big( X^{\xi,\mu}_{t,s}\big) \big], v\big) \mu_s(dv) ds
   \right],
\end{equation}
 where $\mu$ is a relaxed control,
 $\xi$ is an $\mathcal{F}_t$-measurable $\real$-valued random variable,
 $\E_t [X] = \E[X | \mathcal{F}_t]$ is the conditional expectation given
 $\mathcal{F}_t$, $t \in [0,T]$,
 and $(X^{\xi,\mu}_{t,s})_{s\in [t,T]}$ is the
non-linear controlled diffusion given by
\begin{equation}
  \label{main_SDE_relaxed}
\begin{cases}
  \displaystyle
  dX^{\xi,\mu}_{t,s} = \int_U b\big(s,X^{\xi,\mu}_{t,s},v\big) \mu_s(dv) ds + \sigma\big(s,X^{\xi,\mu}_{t,s}\big) dW_s, \quad 0 \leq t<s \leq T,
  \medskip
  \\
X^{\xi,\mu}_{t,t} = \xi.
\end{cases}
\end{equation}
\noindent
 Our approach to the existence of equilibrium controls
 relies on a time discretization of the control problem
using ${n}$-person games, and on
a variation of Pontryagin's maximum principle for
the characterization of
 ${n}$-person games, see Theorem~\ref{MP N person}.
In Corollary~\ref{existence}, we prove the existence of an equilibrium control
in the sense of Definition~\ref{equilibrium_control_relaxed}
for the time-inconsistent mean-field control problem
\eqref{cost_functional_relaxed}-\eqref{main_SDE_relaxed},
 based on
 a formulation of equilibrium controls
as weak limits of the sequence of solutions to
${n}$-person games, see Theorem~\ref{main theorem}.
The proof of Theorem~\ref{main theorem}
uses BSDE convergence arguments and the characterization
Theorem~\ref{MP N person}.

\medskip

The numerical construction of
equilibrium controls
is achieved by approximating ${n}$-person games using finite Markov chains
by adapting the method of \cite{kushner1990numerical} to our setting,
see Theorem~\ref{main theorem approx}.
Precisely, the argument therein 
 applies only to posed $\inf$ problems, as it requires
    comparing
the optimal control $\mu^*$ to any other control $\mu$
 via the inequality $J (t,x,\mu^* ) \leq J(t,x,\mu)$.
Here, the control problem \eqref{cost_functional_relaxed}-\eqref{main_SDE_relaxed}
is not posed inf, instead it is
formulated in the game-theoretic setting of equilibrium controls
in the sense of Definition~\ref{equilibrium_control_relaxed}
below.
 Hence, no such comparison of equilibrium controls
 is possible
 as in \eqref{eq:equilibrium_control_relaxed},
 nevertheless we are able to apply the comparison argument
 to ${n}$-person games since they are posed $\inf$.
 In Section~\ref{sec:numerical},
 the numerical scheme is implemented using a trinomial tree,
 first on a linear-quadratic model which admits an analytic solution,
 and then on a linear-quartic model which does not have analytic solution.

\medskip
The particular case of mean-variance portfolio selection,
where the cost functional
$J(t,\xi,\mu)$ in \eqref{cost_functional_relaxed}
is given by
$$
   J(t,\xi,\mu) =
 - \mathbb{E}_t \big[ X^{\xi,\mu}_{t,T}\big]
 +
 \frac{\gamma}{2}
 \mathbb{E}_t \big[ \big( X^{\xi,\mu}_{t,T} -  \mathbb{E}_t \big[ X^{\xi,\mu}_{t,T}\big]\big)^2 \big]
$$
 where $\gamma >0$
 has been treated in \cite{czichowsky2013time}
using semimartingale theory
for the convergence of equilibrium controls  from discrete to continuous time.
See also \cite{huang2018strong}
 in the case where $\big(X^{\xi,\mu}_{t,s} \big)_{s\in [t,T]}$
 is a finite Markov chains,
 for time-inconsistent control problems
 with infinite horizon.

\medskip

This paper is organized as follows. After stating the necessary
 preliminaries on equilibrium and
relaxed controls, 
 in Section~\ref{sec:characterization_nperson} we present a characterization of ${n}$-person games
using the maximum principle.
In Section~\ref{sec:convergence_nperson} we show the convergence of the solutions of ${n}$-person games to an equilibrium control, and we obtain in turn the existence of an equilibrium control in Corollary~\ref{existence}.
In Section~\ref{sec:markovchain_approx} we
 deal with the convergence properties of the Markov chain approximation
 for the SDE of ${n}$-person games,
 see Theorem~\ref{main theorem approx}.
In Section~\ref{sec:numerical} we present a numerical application of
 the convergence results obtained in Sections~\ref{sec:convergence_nperson} and \ref{sec:markovchain_approx}.
The proofs of the main Theorems~\ref{MP N person},
Corollary~\ref{existence},  Theorems~\ref{main theorem}
 and \ref{main theorem approx}
   rely on technical lemmas presented in appendix. 

\subsection*{Preliminaries} 

Let $T>0$ be a fixed time horizon and $(\Omega,\mathcal{F}, (\mathcal{F}_t)_{t \in [0,T]}, \P )$ be a filtered probability space satisfying the usual conditions,
where $(\mathcal{F}_t)_{t \in [0,T]}$ is the filtration generated by a standard
Brownian motion $(W_t)_{t\in [0,T]}$.
In the sequel, we let $U$ denote a compact subset of $\mathbb{R}$,
and we denote by ${\mathcal B}([0,T]\times U)$ and ${\mathcal B}(U)$
 the Borel $\sigma$-algebra of $[0,T]\times U$ and $U$, respectively.
\begin{definition}
  The space $\Lambda$
 of deterministic relaxed controls is the set of
 nonnegative measures $\lambda$ on $\mathcal{B} ([0,T]\times U) $ such that
 \begin{equation}
   \label{dfjkls}
 \lambda([0,t] \times U) = t, \qquad t \in [0,T].
\end{equation}
\end{definition}
\noindent
 We also denote by $\lambda_t( \cdot )$
 the density such that 
 $\lambda(dt, dv) = \lambda_t(dv)dt$, $t\in [0,T]$,
 whose existence follows from \eqref{dfjkls}.
\begin{definition}\label{relaxed_control_def}
  \begin{enumerate}[i)]
  \item
    The space $\mathcal{U}([0,T])$ of strict controls over $[0,T]$
    is the set of $({\mathcal F}_t)_{t\in [0,T]}$-adapted
 $U$-valued processes.
\item
 The space ${\mathcal R}([0,T])$ of relaxed controls over $[0,T]$
 is the set of 
 $\Lambda$-valued random variables
 $\lambda$ such that $\lambda([0,t] \times B)$ is ${\cal F}_t$-measurable for all $t \in [0,T]$
 and $B\in {\mathcal B}(U)$.
  \end{enumerate}
\end{definition}
\noindent
We now turn to the definition of equilibrium controls
in a game-theoretic setting,
see 
 \cite{ekeland2006being}.
 Given $\mu , \nu \in \mathcal{R}([0,T])$ two relaxed controls
 and $\varepsilon \in (0,T-t]$, we let
 $\nu \otimes_{t,\varepsilon } \mu$ denote
 the local spike variation of $\mu$,
 defined as
\begin{equation}
\nonumber 
( \nu \otimes_{t,\varepsilon } \mu )_s=
\begin{cases}
  \nu_s, & 0 \leq t \leq s \leq t+\varepsilon ,
  \medskip
  \\
  \mu_s, &
  s \in [0,T] \char`\\ [t, t+\varepsilon  ].
\end{cases}
\end{equation}
\noindent
 As in e.g.
\cite{kushner2013numerical},
\cite{buckdahn2011general},
\cite{djehiche2016characterization}
\cite{bahlali2018relaxed},
we assume the following boundedness and smoothness conditions on the
coefficients and cost functions of the problem
\eqref{cost_functional_relaxed}-\eqref{main_SDE_relaxed}.
\begin{assumption}
\begin{enumerate}[i)]
\item The functions $b, \sigma, h, g, \Phi, \Psi$ are uniformly continuous
        and bounded.
  \item The functions
    $b(t,x,u)$,
    $\sigma(t,x)$,
    $\Phi(x)$,
    $\Psi(x)$ are differentiable with respect to $x$ for all $(t, u)\in [0,T]\times U$,
    and their first order (partial) derivatives
 $\partial_x b(t,x,u)$,
    $\partial_x \sigma(t,x)$,
    $\Phi' (x)$,
    $\Psi' (x)$ are differentiable with respect to $x$ for all $(t, u)\in [0,T]\times U$,
    are uniformly continuous and bounded.
  \item The functions $h(t,x,y,u), g(x,y)$ are differentiable with respect to
    $( x, y )$ for all $(t, u)\in [0,T]\times U$,
    and their first order partial derivatives
    $\partial_x h(t,x,y,u)$,
    $\partial_y h(t,x,y,u)$,
    $\partial_x g(x,y)$,
    $\partial_y g(x,y)$
    are uniformly continuous and bounded.
\item There is a constant $\sigma_0 > 0$ such that $\sigma(t,x) \geq \sigma_0$ for all $(t, x)\in [0,T] \times \real$.
\end{enumerate}
\label{basic assumptions}
\end{assumption}
\noindent
We note that the functions $b(t,\cdot ,u)$,
$\sigma(t,\cdot )$,
$h(t,\cdot ,\cdot ,u)$,
$g(\cdot ,\cdot )$,
$\Phi(\cdot )$,
 $\Psi(\cdot )$ are
globally Lipschitz continuous for all $(t, u)\in [0,T]\times U$ 
since they have bounded derivatives. In the sequel, we fix an initial condition
$x_0\in \real$,
and given $\mu \in \mathcal{R}([0,T])$ we
let $X^\mu_t:=X^{x_0 , \mu}_{0,t}$, $t\in [0,T]$,
 denote the solution of the SDE
\begin{equation}\label{main_SDE_relaxed_equilibrium}
  \left\{
  \begin{array}{l}
  \displaystyle
  dX^\mu_t = \int_U b\bigl(t,X^\mu_t,v\bigr) \mu_t(dv) dt + \sigma\bigl(t,X^\mu_t\bigr) dW_t, \qquad 0<t \leq T,
\medskip
\\
X^\mu_0 = x_0.
  \end{array}
  \right.
\end{equation}
     The next definition of equilibrium controls
     is an extension of Definition~2.1 in \cite{hu2012time}
     using the space $\mathcal{R}([0,T])$
     of relaxed controls instead
     of the space $\mathcal{U}([0,T])$ of strict controls.
\begin{definition} \label{equilibrium_control_relaxed}
  We say that a relaxed control $\mu^* \in {\cal R}([0,T])$
  is an equilibrium control for the time-inconsistent mean-field control problem \eqref{cost_functional_relaxed}-\eqref{main_SDE_relaxed} if
 \begin{equation}\label{eq:equilibrium_control_relaxed}
\lim\limits_{h \downarrow 0} \frac{J\big(t,X^{\mu^*}_t,\mu^*\big) - J\big(t,X^{\mu^*}_t,\mu \otimes_{t,h} \mu^*\big)}{h} \leq 0,
\quad
\mu \in \mathcal{R}([0,T]), \ \mbox{a.e.} \ t \in [0,T], \ \mathbb{P}\mbox{-}a.s.,
\end{equation}
 where the equilibrium dynamics
 $\big( X^{\mu^*}_t\big)_{t\in [0,T]}$ is
    the
solution of
 \eqref{main_SDE_relaxed_equilibrium}.
\end{definition}
\noindent
 In the literature,
 Definition~\ref{equilibrium_control_relaxed}
 is usually stated in the space $\mathcal{U}([0,T])$ of strict controls
 instead of using the space $\mathcal{R}([0,T])$
 of relaxed controls, see Definition~\ref{relaxed_control_def}.
The relaxed representation of a strict control $u \in \mathcal{U}([0,T])$
is denoted by
\begin{equation}
    \label{relaxed_control_representation}
    \mu(dt,dv) = \mu_t (dv) dt = \delta_{u_t}(dv) dt,
\end{equation}
 where $\delta_x(dv)$ denotes the Dirac measure at $x\in U$.

\medskip

 As the proof of our existence result Corollary~\ref{existence}
 requires the compactness of the control space,
 we choose to work with the space $\mathcal{R}([0,T])$ of relaxed controls
 because it is compact when endowed with the weak topology.
 Examples of control problems which do not admit
 strict equilibrium controls can be constructed based on
 the non compactness of the space $\mathcal{U}([0,T])$ of strict controls,
 see e.g. the Rademacher function example in
 \S~1 of \cite{valadier1994course}.
     In \cite{hu2012time,hu2017time},
     the existence of equilibrium controls is proved
     without requiring the compactness of the control space,
     however this is for the special case of a linear-quadratic structure
     on the SDE and cost functional.
 
\medskip

 For convenience, we introduce the following notation.
 Given $\mu \in \mathcal{R}([0,T])$ a relaxed control of interest,
 for example $\mu^*$ in Theorem~\ref{MP TC}
 or $\mu^{*n} $ in Theorem~\ref{MP N person}
 below,
 for $\varphi = b, \sigma, h, g$ and
 $\gamma = \Phi$, resp. $\Psi$ when $\varphi = h$, resp. $g$,
 we set the notation
\begin{subequations}
	\begin{empheq}[]{align}
  \label{rc1}
  &
  \partial_x\varphi^{\mu}_{t,s} = \int_U \partial_x\varphi (s, X^{\mu}_s, \mathbb{E}_t [\gamma(X^{\mu}_s)], v) \mu_s(dv),
	\\
  \label{rc2}
  &
\partial_y\varphi^{\mu}_{t,s} =  \gamma'(X^{\mu}_s)\mathbb{E}_t\biggl[\int_U
  \partial_y\varphi (s, X^{\mu}_s, \mathbb{E}_t [\gamma(X^{\mu}_s)], v) \mu_s(dv)\biggr],
        \end{empheq}
\end{subequations}
 where $t \leq s \leq T$ and $X^{\mu}$ is defined in \eqref{main_SDE_relaxed_equilibrium}.
 Next, we now introduce the Hamiltonian function
\begin{equation}
  \label{hamiltonian}
H(t,x,y,\mu,p) = p\int_U b (t,x,v) \mu(dv) - \int_U h(t,x,y, v) \mu(dv),
\end{equation}
where $(t,x,y,p) \in [0,T] \times \mathbb{R}^3$,
and $\mu$ is in the collection $\P (U)$ of all probability measures on $U$.
By abuse of notation, we also denote
\begin{equation}
\nonumber 
H(t,x,y,u,p) = pb (t,x,u) - h(t,x,y,u)
\end{equation}
 when the fourth variable in \eqref{hamiltonian} is $u \in U$.
 The next theorem is a direct extension
 to relaxed controls of the characterization
 of strict equilibrium controls
 proved in Theorem~1 of \cite{djehiche2016characterization}
 using the maximum principle, therefore its proof is omitted.
   \begin{theorem}
   \label{MP TC}
  Let $\mu^* \in {\cal R}([0,T])$ denote a
  relaxed control.
  Consider $\partial_{\boldsymbol{\cdot}} b^{\mu^*}_{t,s}$,
  $\partial_{\boldsymbol{\cdot}} \sigma^{\mu^*}_{t,s}$,
  $\partial_{\boldsymbol{\cdot}} h^{\mu^*}_{t,s}$,
 $\partial_{\boldsymbol{\cdot}} g^{\mu^*}_{t,T}$
 given by 
 \eqref{rc1}-\eqref{rc2},
 and let
     $\big(p^{\mu^* }_{t,s}, q^{\mu^* }_{t,s}\big)_{s\in [t,T]}$
 be the solution of the first order adjoint equation
\begin{equation}\label{TC adjoint}
  \left\{
  \begin{array}{l}
    \displaystyle
    dp^{\mu^*}_{t,s} = -\big(p^{\mu^*}_{t,s}\partial_xb^{\mu^*}_{t,s}  + q^{\mu^*}_{t,s}\partial_x\sigma^{\mu^*}_{t,s} - \partial_xh^{\mu^*}_{t,s} - \partial_yh^{\mu^*}_{t,s}\big) ds + q^{\mu^*}_{t,s} dW_s, \quad 0 \leq t \leq s \leq T,
\medskip
\\
\displaystyle
    p^{\mu^*}_{t,T} = - \partial_xg^{\mu^*}_{t,T} - \partial_yg^{\mu^*}_{t,T}.
  \end{array}
  \right.
\end{equation}
Then, $\mu^*$ is an equilibrium control
 for the problem \eqref{cost_functional_relaxed}-\eqref{main_SDE_relaxed}
 if and only if there exists a pair
 $\big(p_{t,s}^{\mu *},q_{t,s}^{\mu *}\big)_{s\in [t,T]}$
 of $(\mathcal{F}_s)_{s\in [t,T]}$-adapted process
 satisfying \eqref{TC adjoint}, and such that
\begin{equation}
\nonumber
H\bigl(t,X^{\mu^*}_t,\Phi(X^{\mu^*}_t),\nu,p^{\mu^*}_{t,t}\bigr) \leq H\bigl(t,X^{\mu^*}_t,\Phi(X^{\mu^*}_t),\mu^*_t,p^{\mu^*}_{t,t}\bigr), \quad
 \nu \in \P (U), \ \text{a.e. } t \in [0,T], \ \P\mbox{-}a.s.
\end{equation}
 \end{theorem}
 \noindent
 In the sequel,
 $C>0$ represents a generic constant which may change from line to line.
\section{Existence of equilibrium controls}
\subsection{Maximum principle characterization of $n$-person games}
\label{sec:characterization_nperson}
In this section, we consider $n$-person games for the construction of
an equilibrium control later in Section~\ref{sec:convergence_nperson}.
In \cite{yong2012time}, equilibrium HJB equations have been used for the
characterization of equilibria via $n$-person games in control
problems without mean-field terms.
Since the extension of this PDE approach to
the mean-field case may not be straightforward,
we propose instead 
to use the maximum principle
for the construction of equilibrium controls.

\medskip

Given ${n}\geq 1$, we consider the sequence $\{t_k = k T/{n}, \ k = 0, 1, \ldots , {n}\}$ with step size $\Delta_n := T/n$.
Theorem~\ref{MP N person} is a characterization
of the solution of the discretization
 of the time-inconsistent mean-field control problem \eqref{cost_functional_relaxed}-\eqref{main_SDE_relaxed}
 into an ${n}$-person game,
 for use in the proofs of Corollary~\ref{existence} and Theorem~\ref{main theorem}.
\begin{theorem}
  \label{MP N person}
\noindent
Let $n\geq 1$.
Under Assumption~\ref{basic assumptions}, suppose that
 the  ${n}$-person discretized
 time-inconsistent mean-field control problem
\begin{equation}
\label{N-person}
J\bigl(t_k,X^{\mu^{*n} }_{t_k}, \mu^{*n} \bigr) = \inf\limits_{\mu \in \mathcal{R}( [t_k,t_{k+1}] )} J\bigl(t_k,X^{\mu^{*n} }_{t_k}, \mu \otimes_{t_k,\Delta_n} \mu^{*n} \bigr),
\qquad
k=0,1,\ldots , {n}-1,
\end{equation}
admits a solution $\mu^{*n} \in\mathcal{R}([0,T])$
and let
$\big(p^{\mu^{*n} }_{t_k,t}$, $q^{\mu^{*n} }_{t_k,t}\big)_{t\in [t_k,T]}$
  be the solution of the first order adjoint equation
  \begin{equation}
    \label{TC adjoint discrete}
  \left\{
  \begin{array}{ll}
    \displaystyle
    dp^{\mu^{*n} }_{t_k,t} =& -\big(p^{\mu^{*n} }_{t_k,t}\partial_xb^{\mu^{*n} }_{t_k,t}  + q^{\mu^{*n} }_{t_k,t}\partial_x\sigma^{\mu^{*n} }_{t_k,t} - \partial_xh^{\mu^{*n} }_{t_k,t} - \partial_yh^{\mu^{*n} }_{t_k,t}\big) dt + q^{\mu^{*n} }_{t_k,t} dW_t, \quad
    t_k \leq t \leq T,
\medskip
 \\
    \displaystyle
    p^{\mu^{*n} }_{t_k,T} =& - \partial_xg^{\mu^{*n} }_{t_k,T} - \partial_yg^{\mu^{*n} }_{t_k,T},
    \qquad k = 0,1,\ldots , n-1.
  \end{array}
  \right.
\end{equation}
  Then we have
\begin{equation}
\label{MP k H compare N person}
 H\bigl(t,X^{\mu^{*n} }_t,\mathbb{E}_{t_k}[\Phi(X^{\mu^{*n} }_t)],\nu,p^{\mu^{*n} }_{t_k,t}\bigr) \leq H\bigl(t,X^{\mu^{*n} }_t,\mathbb{E}_{t_k}[\Phi(X^{\mu^{*n} }_t)],\mu^{*n} _t,p^{\mu^{*n} }_{t_k,t}\bigr),
\end{equation}
$\nu \in \P (U)$, {\em a.e.} $t \in [t_k,t_{k+1}]$, $\P$-$a.s.$,
 $k=0,1,\ldots , n-1$.
\end{theorem}
\begin{Proof}
 We fix $k \in \{0,1,\ldots , {n}-1\}$ and $t \in [t_k ,t_{k+1})$.
 Given $A\in \mathcal{F}_t$ and
 $\nu \in \P (U)$,
 applying Lemma~\ref{l1.1} below
 to the deviated control
 $\mu_s :=\nu \mathbbm{1}_A + \mu^{*n} _s \mathbbm{1}_{\Omega \setminus A}$,
 $s\in [0,T]$,
we have
\begin{align}
  \label{dgffgd}
&
    J\bigl(t_k ,X^{\mu^{*n} }_{t_k },
    \mu \otimes_{t,\varepsilon} \mu^{*n} \bigr)
    -
    J\bigl(t_k ,X^{\mu^{*n} }_{t_k }, \mu^{*n} \bigr)
      \\
\nonumber
  & = \mathbb{E}_{t_k } \left[
    \int_t^{t+\varepsilon} \bigl(
    H\bigl(s,X^{\mu^{*n} }_s,\mathbb{E}_{t_k}[\Phi(X^{\mu^{*n} }_s)],\mu^{*n} _s,p^{\mu^{*n} }_{t_k,s}
      \bigr) -
  H\bigl(s,X^{\mu^{*n} }_s,\mathbb{E}_{t_k}[\Phi(X^{\mu^{*n} }_s)],
              \mu_s
    ,p^{\mu^{*n} }_{t_k,s}\bigr)
    \bigr)
    ds\right]
 + o(\varepsilon),
\end{align}
as $\varepsilon$ tends to zero.
Since $\mu^{*n} $ is a solution of
 \eqref{N-person},
      the deviation $\mu$ of $\mu^{*n}$ in $\mathcal{R}([0,T])$
      over any time period within $[t_k,t_{k+1}]$ will be sub-optimal.
      Therefore, letting $\varepsilon$ tend to $0$,
      the Lebesgue Differentiation Theorem applied to \eqref{dgffgd} yields
\begin{equation}
   \mathbb{E}_{t_k } \left[
     \mathbbm{1}_A
     \bigl( H\bigl(t,X^{\mu^{*n} }_t,\mathbb{E}_{t_k}[\Phi(X^{\mu^{*n} }_t)],\mu^{*n} _t,p^{\mu^{*n} }_{t_k,t}  \bigr) -
  H\bigl(t,X^{\mu^{*n} }_t,\mathbb{E}_{t_k}[\Phi(X^{\mu^{*n} }_t)],
 \nu
    ,p^{\mu^{*n} }_{t_k,t}\bigr)
 \bigr)
   \right] \geq 0,
\end{equation}
   $a.e.$ $t \in [0,T]$.
Since
$A\in \mathcal{F}_t$ is arbitrary, we conclude to \eqref{MP k H compare N person}.
\end{Proof}
\medskip
\noindent
The next lemma,
which has been used in the proof of Theorem~\ref{MP N person},
yields an expansion of the cost functional
$J\big(t_k ,X^{\mu^{*n} }_{t_k }, \mu \otimes_{t,\varepsilon} \mu^{*n} \big)$
in $\varepsilon$.
For $\mu, \nu \in \mathcal{R}([0,T])$
 and $\varphi = b, \sigma, h, g$ we let
\begin{subequations}
	\begin{empheq}[]{align}
\nonumber 
  &
  \delta \varphi^{\nu,\mu}_{t,s} = \int_U \varphi (s, X^{\mu}_s, \mathbb{E}_t [\gamma(X^{\mu}_s)], v) \nu_s(dv) - \int_U \varphi (s, X^{\mu}_s, \mathbb{E}_t [\gamma(X^{\mu}_s) ], v) \mu_s(dv),
	\\
\nonumber 
  &
  \delta \partial_x\varphi^{\nu,\mu}_{t,s} = \int_U \partial_x\varphi (s, X^{\mu}_s, \mathbb{E}_t [\gamma(X^{\mu}_s)], v) \nu_s(dv) - \int_U \partial_x\varphi (s, X^{\mu}_s, \mathbb{E}_t [\gamma(X^{\mu}_s)], v) \mu_s(dv),
        \end{empheq}
\end{subequations}
$t \leq s \leq T$,
 where $\gamma = \Phi$, resp. $\Psi$ when $\varphi = h$, resp. $g$.
\begin{lemma}
  \label{l1.1}
  Under the assumptions of Theorem~\ref{MP N person},
  fix $t \in [0,T)$ and $k \in \{0, 1, \dots, {n}-1\}$ such that $t_k \leq t < t_{k+1}$,
    and let $\mu \in {\cal R}([0,T])$.
   Then, as $\varepsilon > 0$ tends to zero we have the expansion
\begin{equation}
\nonumber 
J\big(t_k ,X^{\mu^{*n} }_{t_k },
 \mu \otimes_{t,\varepsilon} \mu^{*n} \big) = J\big(t_k ,X^{\mu^{*n} }_{t_k }, \mu^{*n} \big)
+ \mathbb{E}_{t_k } \biggl[
  \int_t^{t+\varepsilon } \big(\delta h^{\mu,\mu^{*n} }_{t_k ,s} - \delta b^{\mu,\mu^{*n} }_{t_k ,s} p^{\mu^{*n} }_{t_k ,s}\big)
    ds\biggr] + o(\varepsilon).
\end{equation}
\end{lemma}
\begin{Proof}
  Let $\big( y^{(\varepsilon)}_s\big)_{s\in [t_k , T]}$ denote the solution of the variational equation
\begin{equation}
\label{variational k}
\begin{cases}
    dy^{(\varepsilon)}_s = \bigl(y^{( \varepsilon )}_s\partial_xb^{\mu^{*n} }_{t_k ,s}  + \delta b^{\mu,\mu^{*n} }_{t_k ,s} \mathbbm{1}_{[t, t+\varepsilon]}(s)\bigr) ds + y^{(\varepsilon )}_s\partial_x \sigma^{\mu^{*n} }_{t_k ,s}  dW_s, \quad t_k  \leq s \leq T,
\medskip
\\
y^{(\varepsilon)}_{t_k } = 0,
\end{cases}
\end{equation}
and, for $ (s, u, \theta) \in [t_k,T] \times U \times [0,1] $, let
$\mu^\varepsilon_s := ( \mu \otimes_{t,\varepsilon} \mu^{*n} )_s$,
$\xi^{(\varepsilon)}_s := X^{\mu^\varepsilon}_s - X^{\mu^{*n}}_s$,
$\eta^{(\varepsilon)}_s := \xi^{(\varepsilon)}_s - y^{(\varepsilon)}_s$,
and use the notation
$$
  \begin{cases}
     \partial_x h_{\theta}(s,u) = \partial_x h\big(s, (1-\theta)X^{\mu^{*n}}_s + \theta X^{\mu^\varepsilon}_s
    , \mathbb{E}_{t_k }\big[(1-\theta) \Phi\big(X^{\mu^{*n}}_s\big) + \theta \Phi\big(X^{\mu^\varepsilon}_s\big)\big], u\big),
     \medskip
     \\
 \partial_xg_\theta
 = \partial_x g\big((1-\theta)X^{\mu^{*n}}_T + \theta X^{\mu^\varepsilon}_T
 , \mathbb{E}_{t_k }\big[(1-\theta) \Psi\big(X^{\mu^{*n}}_T\big) + \theta \Psi\big(X^{\mu^\varepsilon}_T\big)\big]\big),
 \medskip
 \\
 \partial_x\Phi_\theta(s)
 = \partial_x \Phi\big( (1-\theta)X^{\mu^{*n}}_s + \theta X^{\mu^\varepsilon}_s\big),
  \end{cases}
  $$
and similarly for
$\partial_yh_\theta(s,u)$,
$\partial_yg_\theta$,
$\partial_x\Psi_\theta$.
We note that by the flow property
$X^{X^{\mu^{*n} }_{t_k}, \mu^{*n}}_{t_k,s} = X^{\mu^{*n} }_s$ the cost functional
in \eqref{cost_functional_relaxed}
rewrites as
$$
J\big(t_k ,X^\mu_{t_k }, \mu \big)
 = \mathbb{E}_{t_k}\left[
   g\big(X^\mu_T, \mathbb{E}_{t_k} \big[\Psi\big( X^\mu_T\big)\big]
   \big)
   +
   \int_{t_k}^T \int_U h\big(s, X^\mu_s, \mathbb{E}_{t_k}\big[\Phi\big( X^\mu_s\big) \big], v\big) \mu_s (dv) ds
   \right].
$$
 By the fundamental theorem of calculus on $[0,1]$
 we have
\begin{align}
  &J\big(t_k ,X^{\mu^{*n} }_{t_k },\mu^\varepsilon \big) - J\big(t_k ,X^{\mu^{*n}}_{t_k}, \mu^{*n} \big)
  \nonumber
  \\
  &= \mathbb{E}_{t_k }\biggl[
    g\bigl(X^{\mu^\varepsilon}_T, \mathbb{E}_{t_k }\big[\Psi\big(X^{\mu^\varepsilon}_T\big)\big]\bigr) -g\bigl(X^{\mu^{*n}}_T, \mathbb{E}_{t_k }\bigl[\Psi\bigl(X^{\mu^{*n}}_T\bigr)\bigr]\bigr) \nonumber
    \\
    &
    \quad +\int_{t_k }^T \biggl(
    \int_U h\bigl(s, X^{\mu^\varepsilon}_s, \mathbb{E}_{t_k }\big[\Phi(X^{\mu^\varepsilon}_s)\big], v\bigr) \mu^\varepsilon _s(dv)
    - \int_U h\big(s, X^{\mu^{*n}}_s, \mathbb{E}_{t_k }\big[\Phi\big(X^{\mu^{*n}}_s\big)\big], v\big) \mu^{*n} _s(dv)\biggr) ds
    \Bigg] \nonumber\\
  &= \mathbb{E}_{t_k }\biggl[
    \xi^{(\varepsilon)}_T\int_0^1 \partial_xg_{\theta} d\theta
    + \mathbb{E}_{t_k }\bigg[
      \int_0^1 \xi^{(\varepsilon)}_T\partial_x\Psi_{\theta} d\theta \bigg]
    \int_0^1 \partial_yg_{\theta} d\theta
     +
      \int_{t_k }^T
      \delta h^{\mu,\mu^{*n} }_{t_k ,s}\mathbbm{1}_{[t,t+\varepsilon]}(s) ds
     \nonumber\\
    &  \quad
      + \int_{t_k }^T \int_U
      \biggl(\xi^{(\varepsilon)}_s\int_0^1 \partial_x h_{\theta}(s,v)d\theta  +
      \mathbb{E}_{t_k }\bigg[
        \int_0^1 \xi^{(\varepsilon)}_s\partial_x\Phi_{\theta}(s) d\theta
        \bigg]
      \int_0^1 \partial_yh_{\theta}(s,v) d\theta \biggr)  \mu^\varepsilon _s(dv)
      ds
      \Bigg]  \nonumber\\
&= \mathbb{E}_{t_k }\biggl[
      \xi^{(\varepsilon)}_T
    \int_0^1 \big(\partial_xg_{\theta} - \partial_xg_0\big) d\theta
    + (\xi^{(\varepsilon)}_T - y^{(\varepsilon)}_T)\partial_xg^{\mu^{*n} }_{t_k ,T}
    + y^{(\varepsilon)}_T \partial_xg^{\mu^{*n} }_{t_k ,T}
    \nonumber\\
    & \quad +
    \mathbb{E}_{t_k }\bigg[
      \int_0^1 \xi^{(\varepsilon)}_T\partial_x\Psi_{\theta}  d\theta
      \bigg]
    \biggl( \int_0^1 \big(\partial_yg_{\theta} - \partial_yg_0\big) d\theta\biggr)
\nonumber\\
&
\quad +\partial_yg_0
\mathbb{E}_{t_k }\bigg[
  \int_0^1 \xi^{(\varepsilon)}_T\big( \partial_x\Psi_{\theta} - \partial_x\Psi_0 \big) d\theta
  \bigg]
+ \partial_yg_0\mathbb{E}_{t_k }\big[\big(\xi^{(\varepsilon)}_T - y^{(\varepsilon)}_T\big)\partial_x\Psi_0\big]  + y^{(\varepsilon)}_T\partial_yg^{\mu^{*n} }_{t_k ,T}
    \nonumber\\ & \quad
    + \int_{t_k }^T \left( \xi^{(\varepsilon)}_s\int_U \int_0^1 \bigl(\partial_xh_{\theta}(s,v) - \partial_xh_0 (s,v)\bigr) d\theta  \mu^\varepsilon _s(dv)
    \right.
    \nonumber\\
    &
\quad  \quad \quad \quad +  \xi^{(\varepsilon)}_s\biggl(
   \int_U \partial_xh_0 (s,v) \mu^\varepsilon _s(dv) - \int_U \partial_xh_0 (s,v) \mu^{*n} _s(dv)\biggr)
    + \big(\xi^{(\varepsilon)}_s - y^{(\varepsilon)}_s\big)\partial_xh^{\mu^{*n} }_{t_k ,s}  + y^{(\varepsilon)}_s\partial_xh^{\mu^{*n} }_{t_k ,s}
  \nonumber\\
  & \quad  \quad \quad \quad +
  \mathbb{E}_{t_k }\biggl[
    \int_0^1 \xi^{(\varepsilon)}_s\partial_x\Phi_{\theta}(s) d\theta
    \biggr]
  \int_U \int_0^1 (\partial_yh_{\theta}(s,v) - \partial_yh_0 (s,v)) d\theta  \mu^\varepsilon _s(dv)
    \nonumber\\
    & \quad \quad  \quad \quad +
    \mathbb{E}_{t_k }\biggl[
      \int_0^1 \xi^{(\varepsilon)}_s\bigl(\partial_x\Phi_{\theta}(s) - \partial_x\Phi_0(s)\bigr)
      d\theta
      \biggr]
    \int_U \partial_yh_0 (s,v) \mu^\varepsilon _s(dv) \nonumber\\
    & \quad \quad \quad  \quad +  \mathbb{E}_{t_k }\bigl[\xi^{(\varepsilon)}_s\partial_x\Phi_0(s) \bigr] \biggl(
    \int_U \partial_yh_0 (s,v) \mu^\varepsilon _s(dv) - \int_U \partial_yh_0 (s,v) \mu^{*n} _s(dv)\biggr) \nonumber\\
    &   \quad  \quad \quad \quad \left.
        + \mathbb{E}_{t_k }\bigl[\big(\xi^{(\varepsilon)}_s - y^{(\varepsilon)}_s\big)\partial_x\Phi_0(s)\bigr]\int_U \partial_yh_0 (s,v) \mu^{*n} _s(dv)  + y^{(\varepsilon)}_s\partial_yh^{\mu^{*n} }_{t_k ,s}  + \delta h^{\mu,\mu^{*n} }_{t_k ,s}\mathbbm{1}_{[t,t+\varepsilon]}(s)
    \right) ds \Bigg]
  \nonumber
  \\
&= \mathbb{E}_{t_k }\biggl[
    y^{(\varepsilon)}_T \partial_xg^{\mu^{*n} }_{t_k ,T}
    + y^{(\varepsilon)}_T \partial_yg^{\mu^{*n} }_{t_k ,T}
    + \int_{t_k }^T \big( y^{(\varepsilon)}_s \partial_xh^{\mu^{*n} }_{t_k ,s}  + y^{(\varepsilon)}_s \partial_yh^{\mu^{*n} }_{t_k ,s} + \delta h^{\mu,\mu^{*n} }_{t_k ,s}\mathbbm{1}_{[t,t+\varepsilon]}(s)
    \big) ds
    \biggr]
    + o\big(\varepsilon\big), \label{l1.1 last inequality}
  \\
&= \mathbb{E}_{t_k }\biggl[
    \int_{t_k }^T \big( y^{(\varepsilon)}_s \partial_xh^{\mu^{*n} }_{t_k ,s}  + y^{(\varepsilon)}_s \partial_yh^{\mu^{*n} }_{t_k ,s} + \delta h^{\mu,\mu^{*n} }_{t_k ,s}\mathbbm{1}_{[t,t+\varepsilon]}(s)
    \big) ds
    \biggr]
  -
  \mathbb{E}_{t_k }\big[
   y^{(\varepsilon)}_T p^{\mu^{*n} }_{t_k ,T}
   \bigr] + o\big(\varepsilon\big), \label{l1.1 last inequality-1}
\end{align}
 as $\varepsilon$ tends to zero,
 where \eqref{l1.1 last inequality} is due to
 Relations~\eqref{x delta k}, \eqref{x delta y epsilon k} in Lemma~\ref{l1},
 the conditional H\"older inequality, Assumption~\ref{basic assumptions},
 and Lemma~\ref{lemma: uniform and lipschitz},
 and \eqref{l1.1 last inequality-1} is due to \eqref{TC adjoint discrete}.
 We conclude using the identity
$$
  \mathbb{E}_{t_k } \big[y^{(\varepsilon)}_T p^{\mu^{*n} }_{t_k ,T} \big] = \mathbb{E}_{t_k } \left[\int_{t_k }^T \big(y^{(\varepsilon)}_s \partial_xh^{\mu^{*n} }_{t_k ,s} + y^{(\varepsilon)}_s\partial_yh^{\mu^{*n} }_{t_k ,s}
  + \delta b^{\mu,\mu^{*n} }_{t_k ,s}\mathbbm{1}_{[t,t+\varepsilon]}(s) p^{\mu^{*n} }_{t_k ,s}\big) ds\right],
  $$
 that follows from It\^o's lemma.
\end{Proof}
\noindent
In the next lemma, we derive the order of convergence for the variational equation
\eqref{variational k}, which has been used in the proof of Lemma~\ref{l1.1}.
\begin{lemma}
  \label{l1}
  Under the assumptions of Theorem~\ref{MP N person},
  fix $t \in [0,T)$ and $k \in \{0, 1, \dots, {n}-1\}$ such that $t_k \leq t < t_{k+1}$,
    let $\mu \in {\cal R}([0,T])$, $\varepsilon > 0$, $p \geq 1$,
    and denote
 $\xi^{(\varepsilon)}_s = X^{\mu^\varepsilon}_s - X^{\mu^{*n}}_s$, 
$\eta^{(\varepsilon)}_s = \xi^{(\varepsilon)}_s - y^{(\varepsilon)}_s$
    as in \eqref{variational k},
    where
    $\mu^\varepsilon = \mu \otimes_{t,\varepsilon} \mu^{*n}$.
    Then, as $\varepsilon$ tends to zero we have the estimates
\begin{align}
  & \mathbb{E}_{t_k}\biggl[\sup\limits_{s \in [t_k,T]}
    \big| \xi^{(\varepsilon)}_s\big|^{2p}\biggr] = O(\varepsilon^{2p}), \label{x delta k}\\
  &
   \mathbb{E}_{t_k}\biggl[\sup\limits_{s \in [t_k,T]} \abs{y^{(\varepsilon)}_s}^{2p}\biggr] = O(\varepsilon^{2p}), \label{y epsilon k}\\
  &
 \mathbb{E}_{t_k}\biggl[\sup\limits_{s \in [t_k,T]} \abs{\xi^{(\varepsilon)}_s - y^{(\varepsilon)}_s}^{2p}\biggr] = o(\varepsilon^{2p}). \label{x delta y epsilon k}
\end{align}
\end{lemma}
\begin{Proof}
$1)$ Proof of \eqref{x delta k}-\eqref{y epsilon k}.
 Letting
$$
  \left\{
  \begin{array}{l}
  \displaystyle
  \tilde{b}^{(\varepsilon )}_s
  =
  \int_U \big(
  b\big(s, X^{\mu^\varepsilon}_s, v\big)
  -
   b\big(s, X^{\mu^{*n}}_s , v\big) \big) \mu^\varepsilon _s(dv)
=
    \int_U \int_0^1 \partial_xb\big(s, (1-\theta)X^{\mu^{*n}}_s + \theta X^{\mu^\varepsilon}_s, v\big) d\theta \mu^\varepsilon _s(dv)
    ,
\medskip
\\
\displaystyle
\tilde{\sigma}^{(\varepsilon )}_s
= \sigma\big(s, X^{\mu^\varepsilon}_s\big)
-
\sigma\big(s, X^{\mu^{*n}}_s \big)
= \int_0^1 \partial_x\sigma\big(s, (1-\theta)X^{\mu^{*n}}_s + \theta X^{\mu^\varepsilon}_s\big) d\theta,
  \end{array}
  \right.
  $$
  by the fundamental theorem of calculus,
  the process $\big(\xi^{(\varepsilon)}_s\big)_{s\in [t_k,T]}$ satisfies the SDE
$$
\begin{cases}
  d\xi^{(\varepsilon)}_s = \big(\xi^{(\varepsilon)}_s
  \tilde{b}^{(\varepsilon )}_s  + \delta b^{\mu,\mu^{*n} }_{t_k ,s}\mathbbm{1}_{[t,t+\varepsilon]}(s)\big) ds + \xi^{(\varepsilon)}_s\tilde{\sigma}^{(\varepsilon )}_s  dW_s, \qquad t_k \leq s \leq T,
\medskip
\\
\xi^{(\varepsilon)}_{t_k} = 0.
\end{cases}
$$
  Next, by the Burkholder-Davis-Gundy inequality, we have
\begin{eqnarray}
  \lefteqn{
     \!   \! \! \! \! \! \! \! \! \! \! \! \! \!
    \mathbb{E}_{t_k}\biggl[\sup\limits_{s \in [t_k,T]} \big| \xi^{(\varepsilon)}_s\big|^{2p} \biggr]
   = \mathbb{E}_{t_k} \biggl[\sup\limits_{s \in [t_k,T]}
       \biggl| \int_{t_k}^s \big(\xi^{(\varepsilon)}_r \tilde{b}^{(\varepsilon )}_r  + \delta b^{\mu,\mu^{*n} }_{t_k ,r}\mathbbm{1}_{[t,t+\varepsilon]}(r)\big) dr + \int_{t_k}^s \xi^{(\varepsilon)}_r \tilde{\sigma}^{(\varepsilon )}_r  dW_r\biggr|^{2p}
    \biggr]
  }
  \nonumber
  \\
  &\leq & C \mathbb{E}_{t_k}\biggl[
    \int_{t_k}^T\big| \xi^{(\varepsilon)}_s \tilde{b}^{(\varepsilon )}_s \big|^{2p} ds
    + \biggl|
    \int_{t_k}^T \big| \delta b^{\mu,\mu^{*n} }_{t_k ,s} \big|
    \mathbbm{1}_{[t,t+\varepsilon]}(s) ds\biggr|^{2p} + \int_{t_k}^T \big| \xi^{(\varepsilon)}_s
     \tilde{\sigma}^{(\varepsilon )}_s \big|^{2p} ds
    \biggr]
  \nonumber 
  \\
  & \leq & C \mathbb{E}_{t_k}\biggl[
    \int_{t_k}^T \big| \xi^{(\varepsilon)}_s\big|^{2p} ds
    \biggr] + C \varepsilon^{2p}
    \label{l1 third inequality}
  \\
  & \leq & C  \int_{t_k}^T\mathbb{E}_{t_k}\biggl[\sup\limits_{r \in [t_k,s]}
    \big| \xi^{(\varepsilon)}_r\big|^{2p} \biggr] ds  + C \varepsilon^{2p} ,\nonumber
\end{eqnarray}
where 
\eqref{l1 third inequality} is due to the boundedness of $b$, $\sigma$ and their derivatives in Assumption~\ref{basic assumptions}.
The proof of \eqref{x delta k} is completed using Gronwall's inequality,
and \eqref{y epsilon k} can be proved similarly.

\noindent
$2)$ Proof of \eqref{x delta y epsilon k}.
By the fundamental theorem of calculus,
the process $\big(\eta^{(\varepsilon)}_s\big)_{s\in [t_k,T]}$ satisfies the SDE
$$
\begin{cases}
d\eta^{(\varepsilon)}_s = \big(\xi^{(\varepsilon)}_s \tilde{b}^{(\varepsilon )}_s  - y^{(\varepsilon)}_s\partial_xb^{\mu^{*n} }_{t_k ,s}  \big) ds + \big(\xi^{(\varepsilon)}_s \tilde{\sigma}^{(\varepsilon )}_s  - y^{(\varepsilon)}_s\partial_x \sigma^{\mu^{*n} }_{t_k ,s} \big) dW_s, \quad t_k \leq s \leq T,
\medskip
\\
\eta^{(\varepsilon)}_{t_k} = 0.
\end{cases}
$$
 As $\varepsilon$ tends to zero, we have
\begin{align}
  & \mathbb{E}_{t_k}\biggl[\sup\limits_{s \in [t_k,T]}
    \big|\eta^{(\varepsilon)}_s\big|^{2p}
    \biggr]
  = \mathbb{E}_{t_k} \biggl[ \sup\limits_{s \in [t_k,T]}
    \biggl|
    \int_{t_k}^s \biggl(
    \xi^{(\varepsilon)}_r
    \biggl(
    \tilde{b}^{(\varepsilon )}_r - \int_U \partial_xb\bigl(r, X^{\mu^{*n}}_r, v\bigr) \mu^\varepsilon _r(dv)\biggr)
    \nonumber\\
    & \quad
     + \xi^{(\varepsilon)}_r\biggl(
    \int_U \partial_xb\bigl(r, X^{\mu^{*n}}_r, v\bigr) \mu^\varepsilon _r(dv) - \partial_xb^{\mu^{*n} }_{t_k ,r}\biggr)
    + \eta^{(\varepsilon)}_r\partial_xb^{\mu^{*n} }_{t_k ,r} \biggr) dr
    \nonumber\\
    & \left.
    \quad
      + \int_{t_k}^v \Bigl(\xi^{(\varepsilon)}_r\bigl( \tilde{\sigma}^{(\varepsilon )}_r - \partial_x \sigma^{\mu^{*n} }_{t_k ,r}\bigr)
    +
    \eta^{(\varepsilon)}_r\partial_x \sigma^{\mu^{*n} }_{t_k ,r} \Bigr) dW_r\biggr|^{2p}
    \right]
  \nonumber\\
  & \leq C \mathbb{E}_{t_k}
  \left[
    \int_{t_k}^T \biggl(\big|\xi^{(\varepsilon)}_s\big|^{2p}\left|
     \tilde{b}^{(\varepsilon )}_s - \int_U \partial_xb\bigl(s, X^{\mu^{*n}}_s, v\bigr) \mu^\varepsilon _s(dv)
    \right|^{2p} + \bigl|\eta^{(\varepsilon)}_s\partial_xb^{\mu^{*n} }_{t_k ,s}\bigr|^{2p} + \bigl|\eta^{(\varepsilon)}_s\partial_x\sigma^{\mu^{*n} }_{t_k ,s}\bigr|^{2p}
    \right.
    \nonumber\\
    &
    \left.
    \quad
    + \big|\xi^{(\varepsilon)}_s\big|^{2p}\big| \tilde{\sigma}^{(\varepsilon )}_s - \partial_x \sigma^{\mu^{*n} }_{t_k ,s}\big|^{2p} \biggr) ds
    +
    \left|\int_{t_k}^T
\xi^{(\varepsilon)}_s\abs{\delta \partial_xb^{\mu,\mu^{*n} }_{t_k ,s}}\times\mathbbm{1}_{[t,t+\varepsilon]}(s) ds \right|^{2p} \right] \nonumber 
    \\
    & \leq    C
    \biggl(
    \mathbb{E}_{t_k}\biggl[
      \int_{t_k}^{T} \big|\eta^{(\varepsilon )}_s\big|^{2p} ds
      \biggr]
    + \sup\limits_{s \in [t_k,T]}
    \sqrt{\mathbb{E}_{t_k}\big[
        \big| \xi^{(\varepsilon )}_s \big|^{4p} \big]\mathbb{E}_{t_k}\biggl[
        \biggl|
                 \tilde{b}^{(\varepsilon )}_s - \int_U \partial_xb(s, X^{\mu^{*n}}_s, v) \mu^{(\varepsilon )}_s(dv)
            \biggr|^{4p} \biggr] }
    \nonumber
    \\
    &
    \quad
    + \sup\limits_{s \in [t_k,T]}
    \sqrt{ \mathbb{E}_{t_k}\bigl[ \big|
            \xi^{(\varepsilon )}_s
\big|^{4p} \bigr]\mathbb{E}_{t_k}\bigl[ \bigl|
                \tilde{\sigma}^{(\varepsilon )}_s - \partial_x \sigma^{\mu^{*n} }_{t_k ,s}
            \bigr|^{4p} \bigr]}
    + \varepsilon^{4p} \biggr)
    \label{l1(2) third inequality}\\
    & \leq C \int_{t_k}^T \mathbb{E}_{t_k}\biggl[ \sup\limits_{r \in [t_k,s]}
    \big|\eta^{(\varepsilon )}_r\big|^{2p} \biggr] ds  + o\bigl(\varepsilon^{2p}\bigr) ,\label{l1(2) forth inequality}
\end{align}
where  \eqref{l1(2) third inequality} is due to
 \eqref{x delta k} and to the boundedness of
 $\partial_xb$ and $\partial_x\sigma$,
 \eqref{l1(2) forth inequality} is due to the uniform continuity of $\partial_xb$,
  $\partial_x\sigma$ in Assumption~\ref{basic assumptions}, Lemma~\ref{lemma: uniform and lipschitz}, and \eqref{x delta k}.
The proof of \eqref{x delta y epsilon k} is completed by Gronwall's inequality.
\end{Proof}
\subsection{Construction of equilibrium controls}
\label{sec:convergence_nperson}
We equip
the space $\Lambda$ of deterministic relaxed controls with the
 weak topology generated by the bounded continuous functions
on $[0,T]\times U$.
 The spaces $\mathcal{C}([0,T])$ and $\mathcal{D}([0,T])$
 of continuous and c\`adl\`ag functions on $[0,T]$
 are equipped with the uniform and Skorokhod metrics,
 respectively.
In Theorem~\ref{main theorem}, we construct an equilibrium control for
  \eqref{cost_functional_relaxed}-\eqref{main_SDE_relaxed}
  as the weak limit of the solution of the $n$-person game \eqref{N-person}
as ${n}$ tends to infinity, and in Corollary~\ref{existence}
we prove the existence of an equilibrium control.
In the sequel we let
$\floor*{T}_{n} := T$ and $\floor{t}_{n} := t_k$ if $t_k\leq t < t_{k+1}$,
$k=0,1,\ldots , n-1$
\begin{theorem} \label{main theorem}
      Under Assumption~\ref{basic assumptions},
      for any ${n} \geq 1$ there exists a solution $\mu^{*n}$ of
      the $n$-person game \eqref{N-person}
      in the space of relaxed controls.
      In addition,
  the weak limit $\mu^*$ of any convergent subsequence
  of $(\mu^{*n} )_{{n} \geq 1}$ 
  is an equilibrium control for the time-inconsistent mean-field control problem
  \eqref{cost_functional_relaxed}-\eqref{main_SDE_relaxed}.
  \end{theorem}
\begin{Proof}
 We start by constructing a solution $\mu^{*n}$ of
 the $n$-person game \eqref{N-person}
 using backward induction
 in the compact space of relaxed controls. 
 By Theorem~2.14 of \cite{bahlali2018relaxed}
there exists a mapping $\hat{\mu}_n : \real \to \mathcal{R}([t_{n-1}, T])$ such that
\begin{equation*}
    J\bigl(t_{n-1}, x, \hat{\mu}_n (x) \bigr)
    = \inf\limits_{\mu \in \mathcal{R}([t_{n-1},T])}
    J\bigl(t_{n-1}, x, \mu \bigr),
    \qquad x \in \real. 
\end{equation*}
Next, applying this argument recursively to 
$k=n-1,\ldots ,2, 1$, we obtain a mapping
 $\mu_k: \real \to \mathcal{R}([t_{k-1}, t_k])$ such that
\begin{eqnarray*}
  \lefteqn{
    \! \! \! \! \! \! \! \! \! \! \! \! \! 
    J\big(t_{k-1}, x, \mu_k(x)
    \otimes_{t_{k-1},\Delta_n} \hat{\mu}_{k+1} \big(X^{x,\mu_k(x)}_{t_{k-1}, t_k}\big) \big)
  }
  \\
   & & ~~~~~~~~~~~~~  =\inf\limits_{\mu \in \mathcal{R}([t_{k-1}, t_k])}
        J\big(t_{k-1}, x, \mu
        \otimes_{t_{k-1},\Delta_n} \hat{\mu}_{k+1} \big(X^{x,\mu}_{t_{k-1}, t_k}\big) \big),
\end{eqnarray*}
 and let $\hat{\mu}_k(x) := \mu_k(x)
\otimes_{t_{k-1},\Delta_n} \hat{\mu}_{k+1} \big(X^{x,\mu_k(x)}_{t_{k-1}, t_k}\big)$,
 $x \in \real$. 
Then, $\mu^{*n}:=\hat{\mu}_1(x_0) \in \mathcal{R}([0, T])$ is a solution of the n-person game \eqref{N-person}.

  \medskip

  By abuse of notation, we denote by $(\mu^{*n} )_{{n}\geq 1}$ the extracted convergent subsequence on $\Lambda$, and show that its weak limit $\mu^*$ is an equilibrium control.
 We have
\begin{align}
  &\mathbb{E}\bigg[\int_0^T\big| H\bigl(t,X^{\mu^{*n} }_t,\mathbb{E}_{\floor{t}_{n}}\big[\Phi\big(X^{\mu^{*n} }_t\big)\big],\mu^{*n} _t,p^{\mu^{*n} }_{\floor{t}_{n},t}\bigr) - H\bigl(t,X^{\mu^*}_t,\Phi(X^{\mu^*}_t),\mu^*_t,p^{\mu^*}_{t,t}\bigr)\big| dt\bigg]  \nonumber\\
  &= \mathbb{E}\bigg[\int_0^T\left|
    p^{\mu^{*n} }_{\floor{t}_{n},t}\int_U b\big(t,X^{\mu^{*n} }_t,v\big) \mu^{*n} _t(dv)  - \int_U h\big(t,X^{\mu^{*n} }_t,\mathbb{E}_{\floor{t}_{n}}\big[\Phi\big(X^{\mu^{*n} }_t\big)\big], v\big) \mu^{*n} _t(dv)
    \right.
    \nonumber\\
    & \quad
    \left.\left.
    - p^{\mu^*}_{t,t}\int_U b(t,X^{\mu^*}_t,v) \mu^*_t(dv)  + \int_U h(t,X^{\mu^*}_t,\mathbb{E}_t[\Phi(X^{\mu^*}_t)], v) \mu^*_t(dv)\right| dt\right]
  \nonumber\\
  &\leq C \mathbb{E}\bigg[\int_0^T \bigg(
    \bigl|
    p^{\mu^{*n} }_{\floor{t}_{n},t}-p^{\mu^*}_{t,t}\bigr| \times
   \bigg| \int_U b\big(t,X^{\mu^{*n} }_t,v\big) \mu^{*n} _t(dv)\biggr|
 + \big|p^{\mu^*}_{t,t}\big|  \int_U \big| b\big(t,X^{\mu^{*n} }_t,v\big)-b(t,X^{\mu^*}_t,v) \big| \mu^{*n} _t(dv) \nonumber\\
& \quad + \big|p^{\mu^*}_{t,t}\big| \times \bigg| \int_U b(t,X^{\mu^*}_t,v) \mu^{*n} _t(dv)  - \int_U b(t,X^{\mu^*}_t,v) \mu^*_t(dv) \bigg| \nonumber\\
& \quad +  \int_U \big| h\big(t,X^{\mu^{*n} }_t,\mathbb{E}_{\floor{t}_{n}}\big[\Phi\big(X^{\mu^{*n} }_t\big)\big], v\big) - h(t,X^{\mu^*}_t,\mathbb{E}_t[\Phi(X^{\mu^*}_t)], v)| \mu^{*n} _t(dv)  \nonumber\\
    &\left.\left.
     \quad + \bigg| \int_U h(t,X^{\mu^*}_t,\mathbb{E}_t[\Phi(X^{\mu^*}_t)], v) \mu^{*n} _t(dv) - \int_U h(t,X^{\mu^*}_t,\mathbb{E}_t[\Phi(X^{\mu^*}_t)], v) \mu^*_t(dv) \bigg| \right)
    dt\right] \nonumber\\
  &\leq C \int_0^T \mathbb{E} \biggl[\bigg| \int_U h(t,X^{\mu^*}_t,\mathbb{E}_t[\Phi(X^{\mu^*}_t)], v) \mu^{*n} _t(dv) - \int_U h(t,X^{\mu^*}_t,\mathbb{E}_t[\Phi(X^{\mu^*}_t)], v) \mu^*_t(dv) \bigg| \label{main theorem last inequality-0}\\
    & \quad \quad \quad \quad \quad \quad
    + \big|p^{\mu^{*n} }_{\floor{t}_{n},t} - p^{\mu^*}_{t,t}\big| + 2 \big| X^{\mu^{*n} }_t - X^{\mu^*}_t\big|
    + \big|\mathbb{E}_{\floor{t}_{n}}[\Phi(X^{\mu^*}_t)]-\mathbb{E}_t[\Phi(X^{\mu^*}_t)]\big|
    \biggr]dt
  \label{main theorem last inequality-1}
 \\
    &
 \hskip-0.1cm
 + C
  \biggl( \mathbb{E}\left[ \int_0^T\big| p^{\mu^*}_{t,t}\big|^2  dt\right]
  \hskip-0.1cm
  \mathbb{E} \biggl[ \int_0^T \hskip-0.2cm
\bigg(            \big|X^{\mu^{*n} }_t - X^{\mu^*}_t\big|^2 +
            \bigg|
            \int_U b(t,X^{\mu^*}_t,v) \mu^{*n} _t(dv) - \int_U b(t,X^{\mu^*}_t,v) \mu^*_t(dv)\bigg|^2 \bigg)
            dt \bigg] \biggr)^{ 1/2 }.
  \label{main theorem last inequality}
\end{align}
Since
$$\left| \int_U h(t,X^{\mu^*}_t,\mathbb{E}_t[\Phi(X^{\mu^*}_t)], v) \mu^{*n} _t(dv) - \int_U h(t,X^{\mu^*}_t,\mathbb{E}_t[\Phi(X^{\mu^*}_t)], v) \mu^*_t(dv) \right|
$$
is uniformly bounded by some $K>0$
from Assumption~\ref{basic assumptions},
 \eqref{main theorem last inequality-0} converges to $0$ as $n \to \infty$ by
dominated convergence and Lemma~\ref{convergence implies uniform}.
 The first term in \eqref{main theorem last inequality-1} converges to $0$ as $n \to \infty$ by Lemma~\ref{x L2 convergence proof} and
dominated convergence, since by
Theorem~4.2.1 in \cite{zhangjianfeng} there exists $C>0$
  such that
\begin{eqnarray}
\nonumber
  \lefteqn{
      \! \! \! \! \! \! \! \!
      \! \! \! \! \! \! \! \!
      \! \! \! \! \! \! \! \!
      \! \! \! \!
      \sup\limits_{0 \leq t \leq T}  \mathbb{E}\left[
        \sup\limits_{t \leq s \leq T}
        \big(
        \big| p^{{\mu}^n }_{\floor{t}_{n},s}\big|^2 + \big| p^{{\mu}}_{t,s} \big|^2
        \big)
        \right]
\leq C \sup\limits_{0 \leq t \leq T}\mathbb{E}\Biggl[
     \big| \partial_xg^{{\mu}^n }_{\floor{t}_{n},T} + \partial_yg^{{\mu}^n }_{\floor{t}_{n},T}\big|^2 + \big| \partial_xg^{{\mu}}_{t,T} + \partial_yg^{{\mu}}_{t,T} \big|^2
}
\\
  \label{bd}
& &
     + \left(
     \int_t^T \big| \partial_xh^{{\mu}^n }_{\floor{t}_{n},s} + \partial_yh^{{\mu}^n }_{\floor{t}_{n},s} \big| ds\right)^2 + \left(
     \int_t^T \big| \partial_xh^{{\mu}}_{t,s} + \partial_yh^{{\mu}}_{t,s} \big| ds
     \right)^2\Biggr]
\end{eqnarray}
which is bounded uniformly in $t \in [0,T]$ by Assumption~\ref{basic assumptions}.
 The second term in \eqref{main theorem last inequality-1} converges to $0$ as $n \to \infty$ by
Lemma~\ref{x L2 convergence proof}.
The third term in \eqref{main theorem last inequality-1} converges to $0$ as $n \to \infty$ by
Theorem~4 in \cite{fetter1977continuity} and dominated convergence.
 Since
$\big| \int_U b(t,X^{\mu^*}_t,v) \mu^{*n} _t(dv) - \int_U b(t,X^{\mu^*}_t,v) \mu^*_t(dv)\big|^2$ is uniformly bounded by $K^2$, \eqref{main theorem last inequality} converges to $0$ as $n$ tends to infinity by
 Lemma~\ref{x L2 convergence proof}, \eqref{bd} and dominated convergence,
  and Lemma~\ref{convergence implies uniform}.
 Therefore, we have
 \begin{equation}
   \nonumber 
  \lim\limits_{{n}\to \infty}\mathbb{E}\bigg[\int_0^T\big| H\bigl(t,X^{\mu^{*n} }_t,\mathbb{E}_{\floor{t}_{n}}\big[\Phi\big(X^{\mu^{*n} }_t\big)\big],\mu^{*n} _t,p^{\mu^{*n} }_{\floor{t}_{n},t}\bigr) - H\bigl(t,X^{\mu^*}_t,\Phi( X^{\mu^*}_t ),\mu^*_t,p^{\mu^*}_{t,t}\bigr)\big| dt\bigg] = 0,
\end{equation}
 similarly, for any $\nu \in \P (U)$ we have
 \begin{equation}
   \nonumber 
  \lim\limits_{{n}\to \infty}\mathbb{E}\bigg[\int_0^T\big| H\bigl(t,X^{\mu^{*n} }_t,\mathbb{E}_{\floor{t}_{n}}\big[\Phi\big(X^{\mu^{*n} }_t\big)\big],\nu,p^{\mu^{*n} }_{\floor{t}_{n},t}\bigr) - H\bigl(t,X^{\mu^*}_t,\Phi( X^{\mu^*}_t ),\nu,p^{\mu^*}_{t,t}\bigr)\big| dt\bigg] = 0,
\end{equation}
therefore there exists an increasing sequence
 $({n_i})_{{i} \geq 1}$ of integers such that
\begin{equation*} 
  \lim\limits_{{i}\to \infty} H\bigl(t,X^{\mu^{*n_i} }_t,\mathbb{E}_{\floor{t}_{n_i}}\big[\Phi\big(X^{\mu^{*n_i} }_t\big)\big],\mu^{*n_i}_t,p^{\mu^{*n_i}}_{\floor{t}_{{n_i}},t}\bigr) = H\bigl(t,X^{\mu^*}_t,\Phi( X^{\mu^*}_t ),\mu^*_t,p^{\mu^*}_{t,t}\bigr),
\end{equation*}
and
$$
\lim\limits_{{i}\to \infty} H\bigl(t,X^{\mu^{*n_i} }_t,\mathbb{E}_{\floor{t}_{n_i}}\big[\Phi\big(X^{\mu^{*n_i} }_t\big)\big],\nu,p^{\mu^{*n_i}}_{\floor{t}_{{n_i}},t}\bigr)= H\bigl(t,X^{\mu^*}_t,\Phi( X^{\mu^*}_t ),\nu,p^{\mu^*}_{t,t}\bigr),
$$
$a.e.$ $t \in [0,T]$, $\P$-$a.s.$.
 In addition, by Theorem~\ref{MP N person} we have
\[
H\bigl(t,X^{{n_i}}_t,\mathbb{E}_{\floor{t}_{n_i}}\big[\Phi\big(X^{\mu^{*n_i}}_t\big)\big],\nu,p^{\mu^{*n_i}}_{\floor{t}_{n_i},t}\bigr) \leq H\bigl(t,X^{n_i}_t,\mathbb{E}_{\floor{t}_{n_i}}\big[\Phi\big(X^{\mu^{*n_i}}_t\big)\big],\mu^{*n_i}_t,p^{\mu^{*n_i}}_{\floor{t}_{n_i},t}\bigr),
\]
 a.e. $t \in [0,T]$, $\P$-$a.s.$
for all $\nu \in \P (U)$, hence as $k$ tends to infinity we find
\[
H\bigl(t,X^{\mu^{*n_i} }_t,\Phi( X^{\mu^*}_t ),\nu,p^{\mu^*}_{t,t}\bigr) \leq H\bigl(t,X^{\mu^*}_t,\Phi( X^{\mu^*}_t ),\mu^*_t,p^{\mu^*}_{t,t}\bigr), \quad\text{a.e. } t \in [0,T], \ \P\mbox{-}a.s.,
\]
for all $\nu \in \P (U)$,
hence the weak limit $\mu^*$
 of $(\mu^{*n_i})_{i \geq 1}$
 on $\Lambda$
 is an equilibrium control by Theorem~\ref{MP TC}.
\end{Proof}
\noindent
Applying Theorem~2.14 of \cite{bahlali2018relaxed}
under Assumption~\ref{basic assumptions} and using backward induction, for any ${n} \geq 1$
we construct a solution $\mu^{*n}$ of
the $n$-person game \eqref{N-person}
by recursively solving Problem~\eqref{N-person}
in the space $\mathcal{R}( [t_k ,t_{k+1}] )$
of relaxed controls, $k={n}-1, \ldots ,1 , 0$.
By the discussion below
Definition~2.1 in \cite{el1987compactification},
the vague topology used therein on $\Lambda$ 
is equivalent to the weak topology,
and $\Lambda$ is a compact metrizable space
since the set $[0,T]\times U$ is compact.
Therefore, the sequence
 $(\mu^{*n} )_{{n} \geq 1}$ of relaxed controls
 solutions to the 
 $n$-person game \eqref{N-person} is tight, 
 and it admits at least one weakly convergent subsequence, see Theorem~5.1
in \cite{billingsley1999}.
As a consequence, we obtain
the next existence result from Theorem~\ref{main theorem}.
\begin{corollary}
\label{existence}
 Under Assumption~\ref{basic assumptions},
 the time-inconsistent mean-field control problem
 \eqref{cost_functional_relaxed}-\eqref{main_SDE_relaxed}
 admits an equilibrium control $\mu^*$.
\end{corollary}
\begin{Proof}
 By Theorem~\ref{main theorem} above,
 the weak limit $\mu^*$ of any weakly convergent subsequence
 of $(\mu^{*n} )_{{n} \geq 1}$
  is an equilibrium control.
\end{Proof}
\noindent
\noindent
Applying Theorem~\ref{main theorem} requires to check the
weak convergence of a subsequence of $(\mu^{*n} )_{{n}\geq 1}$
in $\Lambda$.
The next corollary shows that this may not be necessary
 if only the value function is concerned.
\begin{corollary}
  \label{main corollary}
  Under Assumption~\ref{basic assumptions},
  the sequence $(J(0, x_0 , \mu^{*n}))_{n \ge 1}$
  admits at least one convergent subsequence.
  In addition, the limit of any such subsequence
  can be written as $J(0, x_0 , \mu^*)$.
\end{corollary}
\begin{Proof}
  Denoting by $( \mu^{*n_i})_{ i \geq 1 }$
  the weakly convergent subsequence of
  $( \mu^{*n})_{ n \geq 1 }$, it suffices to note that
  by Lemmas~\ref{x L2 convergence proof} and \ref{convergence implies uniform},
  the sequence
  $\left( J(0, x_0, \mu^{*n_i})\right)_{ i \geq 1 } $ converges to $J(0, x_0, \mu^*)$
  due to the Lipschitz continuity of $h$ and $g$ in Assumption~\ref{basic assumptions}.
\end{Proof}
\noindent
The next lemma contains stability results for the SDE
\eqref{main_SDE_relaxed_equilibrium}
and for the backward SDE \eqref{TC adjoint discrete}, which have been used in
the proofs of Theorem~\ref{main theorem} and Corollary~\ref{main corollary}.
\begin{lemma}
\label{x L2 convergence proof}
 Let $( {\mu}^n )_{{n}\geq 1} \subset {\cal R}([0,T])$ be a sequence of
  $\Lambda$-valued relaxed controls converging
   weakly to ${\mu} \in {\cal R}([0,T])$.
  Then under Assumption~\ref{basic assumptions},
 we have
\begin{equation}
\nonumber 
\lim\limits_{{n}\to \infty} \mathbb{E}
\left[ \sup\limits_{0 \leq t \leq T} \big| X^{{\mu}^n }_t - X^{{\mu}}_t \big|^2
  \right] = 0 \ \mbox{~and~} \
\lim\limits_{{n}\to \infty} \mathbb{E}\left[ \sup\limits_{t \leq s \leq T}
  \big| p^{{\mu}^n }_{\floor{t}_{n},s} - p^{{\mu}}_{t,s}\big|^2\right] = 0,
\quad t \in [0,T).
\end{equation}
\end{lemma}
\begin{Proof}
 \noindent
 $(i)$
 Using
Assumption~\ref{basic assumptions}, we have
\begin{align*}
  &\mathbb{E}\left[ \sup\limits_{0 \leq t \leq T} \big|X^{{\mu}^n }_t - X^{{\mu}}_t\big|^2\right]
  \\
  &= \mathbb{E}\left[ \sup\limits_{0 \leq t \leq T}
    \left| \int_0^t \int_U b\big(s,X^{{\mu}^n }_s,v\big) {\mu}^n_s(dv) ds - \int_0^t \int_U b(s,X^{{\mu}}_s,v) {\mu}_s(dv) ds
      \right. \right.
      \\
      &\left.\left.
\quad       + \int_0^t \big( \sigma\big(s,X^{{\mu}^n }_s\big) - \sigma(s,X^{{\mu}}_s) \big) dW_s\right|^2\right] \\
  &\leq C \mathbb{E} \left[
    \sup\limits_{0 \leq t \leq T}\bigg(\bigg| \int_0^t \int_U \bigl( b\big(s,X^{{\mu}^n }_s,v\big) - b(s,X^{{\mu}}_s,v) \bigr) {\mu}^n _s(dv) ds\bigg|^2
    \right.
    \\
    &\left. \left.
    \quad + \left| \int_0^t \int_U b(s,X^{{\mu}}_s,v) {\mu}^n_s(dv) ds - \int_0^t \int_U b(s,X^{{\mu}}_s,v) {\mu}_s(dv) ds\right|^2\right)
    + \int_0^T \big| \sigma\big(s,X^{{\mu}^n }_s\big)- \sigma(s,X^{{\mu}}_s)\big|^2 ds\right] \\
  &\leq
  C \mathbb{E}\left[ \int_0^T |f^n _s|^2 ds \right]
   + C \int_0^T \mathbb{E}\left[ \sup\limits_{0 \leq s \leq t}\big|X^{{\mu}^n }_s - X^{{\mu}}_s\big|^2\right] dt,
\end{align*}
where 
$$
f^n _s := \int_U b(s,X^{{\mu}}_s,v) {\mu}^n _s(dv)  - \int_U b(s,X^{{\mu}}_s,v) {\mu}_s(dv).
$$
Hence, by Gronwall's inequality we get
\[
\mathbb{E}\left[ \sup\limits_{0 \leq t \leq T} \big|X^{{\mu}^n }_t - X^{{\mu}}_t\big|^2\right] \leq
 C \re^{CT} \mathbb{E}\left[ \int_0^T |f^n _s|^2 ds \right].
\]
By Lemma~\ref{convergence implies uniform} we have
$\lim_{{n}\to \infty} |f^n _s|^2 = 0$, a.e. $s \in [0,T]$, $\P$-$a.s$,
hence we conclude by dominated convergence as $|f^n _s|^2$ is uniformly bounded by $K^2$.

\noindent
$(ii)$
We fix $t \in [0,T]$ and denote
$$
\left\{
\begin{array}{l}
  \xi^n  := - \partial_xg^{{\mu}^n }_{\floor{t}_{n},T} - \partial_yg^{{\mu}^n }_{\floor{t}_{n},T},
  \medskip
  \\
  \xi := - \partial_xg^{{\mu}}_{t,T} - \partial_yg^{{\mu}}_{t,T},
  \medskip
  \\
  f^n (s,p,q) = p\partial_xb^{{\mu}^n }_{\floor{t}_{n},s}  + q\partial_x \sigma^{{\mu}^n }_{\floor{t}_{n},s}  - \partial_xh^{{\mu}^n }_{\floor{t}_{n},s} - \partial_yh^{{\mu}^n }_{\floor{t}_{n},s},
  \medskip
  \\
f(s,p,q) := p\partial_xb^{{\mu}}_{t,s}  + q\partial_x \sigma^{{\mu}}_{t,s}  - \partial_xh^{{\mu}}_{t,s} - \partial_yh^{{\mu}}_{t,s}.
\end{array}
\right.
$$
 By Theorem~4.4.3 in \cite{zhangjianfeng}, it suffices to show that
 $\lim\limits_{{n}\to \infty} \mathbb{E}\bigl[ |\xi^n  - \xi |^2 \bigr] = 0$,
 that
\begin{equation}
  \label{fjkldsfs1}
  \lim\limits_{{n}\to \infty} \mathbb{E}\biggl[
  \int_t^T | f^n (s,0,0) - f(s,0,0) |^2 ds \biggr] = 0,
\end{equation}
and that $f^n (s,p,q) - f(s,p,q)$ converges to $0$ in $ds \times dP$-measure
as ${n}$ tends to infinity for any fixed $(p,q)$.
We note that the latter condition follows from
 Chebyshev's inequality
 and
\begin{equation}
  \lim\limits_{{n}\to \infty} \mathbb{E}\left[
  \int_t^T
   |f^n (s,p,q) - f(s,p,q)|^2
  ds \right] = 0,
  \quad
  (p,q)\in \real^2.
\end{equation}
Since the arguments leading to the above conditions are similar, we focus
on the limit \eqref{fjkldsfs1}.
 By \eqref{rc2}, we have
\begin{align}
  &\mathbb{E}\left[
    \int_t^T \big| \partial_yh^{{\mu}}_{t,s} - \partial_yh^{{\mu}^n }_{\floor{t}_{n},s}
  \big|^2 ds\right]
  \nonumber\\
  &= \int_t^T \mathbb{E}\Biggl[
    \biggl|
    \Phi'(X^{{\mu}}_s)\mathbb{E}_t \biggl[ \int_U \partial_yh(s, X^{{\mu}}_s, \mathbb{E}_t [\Phi(X^{{\mu}}_s)], v)  {\mu}_s(dv)\biggr]
      \nonumber\\
      & \quad
      - \Phi'\big(X^{{\mu}^n }_s\big)\mathbb{E}_{\floor{t}_{n}} \biggl[
        \int_U \partial_yh\big(s, X^{{\mu}^n }_s, \mathbb{E}_{\floor{t}_{n}} \big[\Phi\big(X^{{\mu}^n }_s\big)\big], v) {\mu}^n _s(dv)\biggr]
  \biggr|^2\Biggr] ds \nonumber\\
  &\leq C \int_t^T \mathbb{E} \Biggl[
    \big| \Phi'(X^{{\mu}}_s) - \Phi'\big(X^{{\mu}^n }_s\big) \big|^2 \times \biggl|\mathbb{E}_t \biggl[
        \int_U \partial_yh(s, X^{{\mu}}_s, \mathbb{E}_t [\Phi(X^{{\mu}}_s)], v) {\mu}_s(dv)
        \biggr]\biggr|^2  +\big| \Phi'\big(X^{{\mu}^n }_s\big) \big|^2
    \nonumber\\
    &\quad\times \biggl( \biggl| \mathbb{E}_t \biggl[ \int_U \partial_yh(s, X^{{\mu}}_s, \mathbb{E}_t [\Phi(X^{{\mu}}_s)], v) {\mu}_s(dv)\biggr]
     - \mathbb{E}_{\floor{t}_{n}} \biggl[
        \int_U \partial_yh(s, X^{{\mu}}_s, \mathbb{E}_t [\Phi(X^{{\mu}}_s)], v) {\mu}_s(dv)\biggr]
      \biggr|^2 \nonumber\\
      &\quad
      + \biggl|
      \mathbb{E}_{\floor{t}_{n}} \biggl[ \int_U \partial_yh(s, X^{{\mu}}_s, \mathbb{E}_t [\Phi(X^{{\mu}}_s)], v) {\mu}_s(dv) - \int_U \partial_yh(s, X^{{\mu}}_s, \mathbb{E}_t [\Phi(X^{{\mu}}_s)], v) {\mu}^n _s(dv)\biggr] \biggr|^2         \nonumber\\
        &\quad
        + \biggl| \mathbb{E}_{\floor{t}_{n}} \biggl[ \int_U \bigl( \partial_yh(s, X^{{\mu}}_s, \mathbb{E}_t [\Phi(X^{{\mu}}_s)], v) - \partial_yh(s, X^{{\mu}}_s, \mathbb{E}_{\floor{t}_{n}} [\Phi(X^{{\mu}}_s)], v)  \bigr) {\mu}^n _s(dv)\biggr]  \biggr|^2
        \nonumber\\
        &\quad
        +\biggl| \mathbb{E}_{\floor{t}_{n}} \biggl[ \int_U \bigl( \partial_yh(s, X^{{\mu}}_s, \mathbb{E}_{\floor{t}_{n}} [\Phi(X^{{\mu}}_s)], v) - \partial_yh\big(s, X^{{\mu}^n }_s, \mathbb{E}_{\floor{t}_{n}} \big[\Phi\big(X^{{\mu}^n }_s\big)\big], v\big)  \bigr) {\mu}^n _s(dv)\biggr] \biggr|^2
        \biggr) \biggr] ds
  \nonumber\\
&\leq C \int_t^T \mathbb{E} \biggl[
    \big| \Phi'(X^{{\mu}}_s) - \Phi'\big(X^{{\mu}^n }_s\big) \big|^2
          \label{p q L2 convergence proof third inequality1}
          \\
    &\quad+  \biggl| \mathbb{E}_t \biggl[ \int_U \partial_yh(s, X^{{\mu}}_s, \mathbb{E}_t [\Phi(X^{{\mu}}_s)], v) {\mu}_s(dv)\biggr]
     - \mathbb{E}_{\floor{t}_{n}} \biggl[
        \int_U \partial_yh(s, X^{{\mu}}_s, \mathbb{E}_t [\Phi(X^{{\mu}}_s)], v) {\mu}_s(dv)\biggr]
     \biggr|^2
     \label{p q L2 convergence proof third inequality2}
     \\
      &\quad
      + \biggl|
      \int_U \partial_yh(s, X^{{\mu}}_s, \mathbb{E}_t [\Phi(X^{{\mu}}_s)], v) {\mu}_s(dv) - \int_U \partial_yh(s, X^{{\mu}}_s, \mathbb{E}_t [\Phi(X^{{\mu}}_s)], v) {\mu}^n _s(dv) \biggr|^2
      \label{p q L2 convergence proof third inequality3}
      \\
        &\quad
        + \biggl| \int_U \bigl( \partial_yh(s, X^{{\mu}}_s, \mathbb{E}_t [\Phi(X^{{\mu}}_s)], v) - \partial_yh(s, X^{{\mu}}_s, \mathbb{E}_{\floor{t}_{n}} [\Phi(X^{{\mu}}_s)], v)  \bigr) {\mu}^n _s(dv)\biggr|^2
        \label{p q L2 convergence proof third inequality4}
        \\
        &\quad
        +\biggl| \int_U \bigl( \partial_yh(s, X^{{\mu}}_s, \mathbb{E}_{\floor{t}_{n}} [\Phi(X^{{\mu}}_s)], v) - \partial_yh\big(s, X^{{\mu}^n }_s, \mathbb{E}_{\floor{t}_{n}} \big[\Phi\big(X^{{\mu}^n }_s\big)\big], v\big)  \bigr) {\mu}^n _s(dv)\biggr|^2
        \Biggr] ds.
  \label{p q L2 convergence proof third inequality}
\end{align}
The inequality \eqref{p q L2 convergence proof third inequality} is due to Assumption~\ref{basic assumptions}, the conditional Jensen's inequality.
Fix $s \in [t,T]$.
By Lemma~\ref{x L2 convergence proof} and Theorem~4 in \cite{fetter1977continuity},
\begin{equation*}
    \lim_{{n}\to \infty} \E\big[
    \big| X^{{\mu}}_s - X^{{\mu}^n }_s\big|^2 + \abs{\mathbb{E}_t [\Phi(X^{{\mu}}_s)] - \mathbb{E}_{\floor{t}_{n}} [\Phi(X^{{\mu}}_s)]}^2 \big] = 0,
\end{equation*}
and
\eqref{p q L2 convergence proof third inequality1},
\eqref{p q L2 convergence proof third inequality4},
\eqref{p q L2 convergence proof third inequality}
converge to zero
by Lemma~\ref{lemma: uniform and lipschitz}, conditional Jensen's inequality, Assumption~\ref{basic assumptions} and dominated convergence on $[0,T]$.
Similarly,
\eqref{p q L2 convergence proof third inequality2} and
\eqref{p q L2 convergence proof third inequality3}
tend to zero by Theorem~4 in \cite{fetter1977continuity} and  Lemma~\ref{convergence implies uniform} respectively.
The term in $\partial_xh^{{\mu}^n }_{\floor{t}_{n},s} - \partial_xh^{{\mu}}_{t,s}$
is treated similarly using \eqref{rc1}.
\end{Proof}
\noindent
\noindent
{\em Remark.}
We note that the equilibrium control $\mu^*$
constructed in Corollary~\ref{existence}
using equal partitions may not be unique.
Indeed, two sequences $(\Pi_1^n)_{n\geq 1}$ and $(\Pi_2^n)_{n\geq 1}$
of partitions may yield distinct limiting equilibrium
controls $\mu_1$ and $\mu_2$
by Theorem~\ref{main theorem}.
 However, under Lipschitz conditions on
the function
\begin{equation*}
    \psi(t,x,p) = {\rm argmax~} H(t,x,\Phi(x),\cdot,p), \quad t\in [0,T], \quad (x,p) \in \real^2,
\end{equation*}
and on the coefficient derivatives
appearing in Assumption~\ref{basic assumptions},
it can be shown by a contraction argument in small time $T$ that
the equilibrium control of
\eqref{cost_functional_relaxed}-\eqref{main_SDE_relaxed}
 can be represented as in \eqref{relaxed_control_representation}
from a strict control in $\mathcal{U}([0,T])$
which is unique in $L^1(\Omega \times [0,T])$.
\section{Numerical implementation}
\subsection{Markov chain approximation of $n$-person games}
\label{sec:markovchain_approx}
 Using Markov chains as in \cite{kushner1990numerical},
 we construct an approximation for
 the relaxed control solution
 $\mu^{*n}$ of the $n$-person game
 \eqref{N-person}  used in Theorem~\ref{main theorem}.
 Then, in Theorem~\ref{main theorem approx} we
 show the convergence of this approximation to
 $\mu^{*n}$, $n\geq 1$.
 Let $n, m \geq 1$, $\Delta_{{n},{m}} := {T}/({n}{m})$,
 and $t^m_k = k \Delta_{n,m}$, $k = 0, 1, \ldots , {n}{m}$.
\begin{definition}
  For any $n, m \geq 1$, we let $\mathcal{U}^{{n},{m}}([0,T])$ denote the set of
  admissible discrete-time strict control sequences $(u_k)_{0 \leq k < nm}$
  such that $u_k$ is ${\cal F}_{t^m_k}$-measurable. 
\end{definition}
\noindent
 Given a sequence $(x_k)_{k = 0,1, \ldots , {n}{m}}$, we let
 $\bar{x}$ be the step function defined as
 $$
 \bar{x}_t = \sum_{k=0}^{nm-1} x_k \mathbbm{1}_{[t^m_k,t^m_{k+1})} (t)
   + x^{nm} \mathbbm{1}_{\{ T\}} (t),
 \qquad t\in [0,T].
 $$
 We also let
 $H_{{n},{m}}:\mathbb{R} \to \sqrt{\Delta}_{{n},{m}}\mathbb{Z}$ denote
 the rounding function on 
 $\sqrt{\Delta}_{{n},{m}}\mathbb{Z}$,
      where $\mathbb{Z} = \{ \ldots , -2, -1 , 0 , 1 , 2 , \ldots \}$
      is the set of integers.
\begin{assumption}
\label{discrete_cond}
 Let $n, m \geq 1$, and $u^{{n},{m}} = (u_k)_{0 \leq k < nm}
 \in \mathcal{U}^{{n},{m}}([0,T])$
 be a sequence of admissible discrete-time strict controls.
 We assume that there exists a
discrete-time Markov chain
 $\big( X^{{n},{m},u}_k \big)_{k=0,\ldots, {n}{m}}$
on $\sqrt{\Delta}_{{n},{m}}\mathbb{Z}$,
such that
\begin{enumerate}[(i)]
\item $X^{{n},{m},u}_0 = H_{{n},{m}}(x_0)$,
\item $\P \big( X^{{n},{m},u}_{k+1}  = y \ \big|
  \
  (X^{{n},{m},u}_l, u^{{n},{m}}_l)_{l=0,1,\ldots ,k}
  \big)
  =
  \P \big( X^{{n},{m},u}_{k+1} = y \ \big|
  \
  X^{{n},{m},u}_k, \ u^{{n},{m}}_k \big)$,
  \ $y \in \sqrt{\Delta}_{{n},{m}}\mathbb{Z}$,
\item
  \label{piii}
  $\mathbb{E} \big[ X^{{n},{m},u}_{k+1} - X^{{n},{m},u}_k \ \big|
  \
  (X^{{n},{m},u}_l, u^{{n},{m}}_l)_{l=0,1,\ldots ,k}
  \big] =\Delta_{{n},{m}} b\big(t^m_k , X^{{n},{m},u}_k, u_k\big)$,
  \item
    $\mathbb{E} \big[ \big(X^{{n},{m},u}_{k+1} - X^{{n},{m},u}_k- \Delta_{{n},{m}} b\big(t^m_k, X^{{n},{m},u}_k, u_k\big)\big)^2 \ \big|
      \
  (X^{{n},{m},u}_l, u^{{n},{m}}_l)_{l=0,1,\ldots ,k}
       \big]
    $
    \\
    $=\Delta_{{n},{m}} \sigma^2\big(t^m_k, X^{{n},{m},u}_k\big)+o(\Delta_{{n},{m}})$,
$ k  = 0, 1, \ldots , {n}{m}-1$.
\item There exists $C > 0$ such that
  $\sup\limits_{0 \leq k < nm} \big\vert X^{{n},{m},u}_{k+1} - X^{{n},{m},u}_k \big\vert \leq
  C \sqrt{\Delta}_{{n},{m}}$, ${{n},{m}}\geq 1$.
\end{enumerate}
\end{assumption}
\noindent
 Let
 $\big(\mathcal{F}^{{n},{m}}_t\big)_{t\in [0,T]}$
 denote the filtration generated by $\big(\widebar{X}^{{n},{m},u}_t\big)_{t\in [0,T]}$.
 Given
 $u^{{n},{m}} = (u_k)_{0\leq k < nm} \in \mathcal{U}^{{n},{m}}([0,T])$
 an admissible control sequence,
  we define the cost functional
 \begin{eqnarray}
\nonumber 
  J^{{n},{m}}\bigl(t^m_k,\widebar{X}^{{n},{m},u}_{t^m_k},\bar{u}^{{n},{m}}\bigr) & = &
  \mathbb{E}\bigg[
    g\big(\widebar{X}^{{n},{m},u}_{T}, \mathbb{E}\big[\Psi\big(\widebar{X}^{{n},{m},u}_{T}\big)
      \ \big| \ \mathcal{F}^{{n},{m}}_{t^m_k}\big]\big)
    \\
    \nonumber
    & & \quad +
\int_{t^m_k}^T h\big(s, \widebar{X}^{{n},{m},u}_s, \mathbb{E}\big[\Phi(\widebar{X}^{{n},{m},u}_s)
     \ \big| \ \mathcal{F}^{{n},{m}}_{t^m_k}\big], \bar{u}^{{n},{m}}_s\big) ds
    \ \Big| \ \mathcal{F}^{{n},{m}}_{t^m_k}\bigg],
  \end{eqnarray}
$k=0,1,\ldots , {n}{m}-1$.
Consider the discretization
\begin{equation}
\label{problem NM-person}
 J^{{n},{m}}\big(t_k,\widebar{X}^{{n},{m},u^{*}}_{t_k},\bar{u}^{*{n},{m}}\big) = \inf\limits_{u \in \mathcal{U}^{{n},{m}}( [0,T] )} J^{{n},{m}}\big(t_k,\widebar{X}^{{n},{m},u^{*}}_{t_k}, \bar{u} \otimes_{t_k,\Delta_n } \bar{u}^{*{n},{m}}\big),
\end{equation}
 $k=0,1,\ldots , {n}-1$,
of the ${n}$-person game \eqref{N-person},
which admits a solution $\bar{u}^{*{n},{m}}$ due to the compactness of $U$.
\noindent
 Let the sequence $(W^{{n},{m},u}_k)_{k=0,\ldots, {n}{m}}$ be defined
 by $W^{{n},{m},u}_0 := 0$ and
 \begin{align}
   \nonumber 
  W^{{n},{m},u}_{k+1} - W^{{n},{m},u}_k & : =\frac{X^{{n},{m},u}_{k+1} - X^{{n},{m},u}_k - \Delta_{{n},{m}} b\big(t^m_k, X^{{n},{m},u}_k, u_k\big)}{\sigma\big(t^m_k, X^{{n},{m},u}_k\big)},
  \quad k=0,1,\ldots ,nm-1.
\end{align}
\noindent
 By Assumption~\ref{discrete_cond}-\eqref{piii} we check that
$\big(\widebar{W}^{{n},{m},u}_t\big)_{t\in [0,T]}$ is a martingale with respect to
its own filtration,
 which coincides with
 $\big(\mathcal{F}^{{n},{m}}_t\big)_{t\in [0,T]}$.
 By the Skorokhod representation Theorem~\ref{skorohod theorem},
 all processes can be defined on a same probability space
 $(\Omega, \mathcal{F}, \P )$.
  The next lemma
follows from
 Theorem~4.6 in \cite{kushner1990numerical},
 see also Theorem~10.4.1 in \cite{kushner2013numerical}.
 \noindent
\begin{lemma}\label{Markov chain convergence}
Under Assumptions~\ref{basic assumptions} and \ref{discrete_cond},
  fix $ {n} \geq 1 $ and for any $m\geq 1$
  let
  $u^{{n},{m}} = (u_k)_{0\leq k < nm} \in \mathcal{U}^{{n},{m}}([0,T])$
  be an admissible control sequence.
  Then, letting $\mu^{{n},{m}}$ denote
  the relaxed control representation
  of $\bar{u}^{{n},{m}}$, $m\geq 1$,
  see \eqref{relaxed_control_representation},
  \begin{enumerate}[a)]
    \item
    the sequence
    $\big( \widebar{X}^{{n},{m},u},\mu^{{n},{m}}, \widebar{W}^{{n},{m},u} \big)_{{m}\geq 1}$ is tight
    on $\mathcal{D}([0,T]) \times \Lambda \times \mathcal{D}([0,T])$,
  \item
    the limit of any weakly converging subsequence
    of $\big( \widebar{X}^{{n},{m},u},\mu^{{n},{m}}, \widebar{W}^{{n},{m},u} \big)_{{m}\geq 1}$
takes the form
$(X^\mu ,\mu, W)$ on $\mathcal{D}([0,T]) \times \Lambda \times \mathcal{D}([0,T])$,
where  $W=(W_t)_{t\in [0,T]}$ is a Wiener process and $X^\mu $
  solves \eqref{main_SDE_relaxed_equilibrium} with the relaxed control $\mu$.
  \end{enumerate}
  \end{lemma}
\noindent
The following approximation lemma, 
see e.g. Theorems~3.2.2 and 3.5.2 in \cite{kushner1990weak}
and references therein,
 will be used to approximate relaxed controls $\mu \in {\cal R}([0,T])$
 using elements
 an admissible control sequences in $\mathcal{U}^{{n},{m}}([0,T])$.
\begin{lemma}\label{chattering lemma}
  [Chattering lemma]
  Let $ n \geq 1 $ and $\mu \in \mathcal{R}([0,T])$.
  Under Assumptions~\ref{basic assumptions} and \ref{discrete_cond}
    there exists a sequence
  $(u^{n,m})_{m\geq 1}$
  of admissible controls $u^{n,m}\in \mathcal{U}^{{n},{m}}([0,T])$,
  $m\geq 1$, such that
  the relaxed control representation
  $(\mu^{{n},{m}})_{m \geq 1}$ of $(\bar{u}^{{n},{m}})_{m\geq 1}$,
  see \eqref{relaxed_control_representation},
  converges weakly to $\mu$ on $\Lambda$ as $m$ tends to infinity.
\end{lemma}
\begin{Proof}
   The sequence $(\bar{u}^{{n},{m}})_{m\geq 1}$ is constructed in the proof of
  Theorem~3.5.2 in \cite{kushner1990weak}
  and its relaxed control representation
  $(\mu^{{n},{m}})_{m \geq 1}$
  is shown to converge weakly to $\mu$.
\end{Proof}
\noindent
The next theorem, which is the main result of this section,
shows the convergence of the solution of
the discretized problem \eqref{problem NM-person}
  to the solution of the $n$-person game \eqref{N-person}.
\begin{theorem}
  \label{main theorem approx}
  Under Assumption~\ref{basic assumptions},
  fix $ {n} \geq 1 $ and
  let $(\bar{u}^{*{n},{m}})_{{m} \geq 1}$ be a sequence of solutions
  to Problem~\eqref{problem NM-person}
  with relaxed control representation
  $(\mu^{*{n},{m}})_{{m} \geq 1}$, see \eqref{relaxed_control_representation},
  and let $(\widebar{X}^{{n},{m},u^*})_{{m} \geq 1}$
  denote the Markov chain
  defined in Assumption~\ref{discrete_cond}. 
 Then,
  \begin{enumerate}[a)]
      \item The sequence $\big(\widebar{X}^{{n},{m},u^*},\mu^{*{n},{m}} \big)_{{m}\geq 1}$ is tight on $\mathcal{D}([0,T]) \times \Lambda$.
\item
  \label{b}
  Denoting by $(X^{\mu^{*n}},\mu^{*n})$ the limit of
  any weakly converging subsequence of
  $\big( \widebar{X}^{{n},{m},u^*},\mu^{*{n},{m}} \big)_{{m}\geq 1}$ on $\mathcal{D}([0,T]) \times \Lambda$,
  the process $(X^{\mu^{*n}}_t)_{t\in [0,T]}$ solves the SDE
  \eqref{main_SDE_relaxed_equilibrium} with relaxed control $\mu^{*n}$.
\item
  \label{c}
  The relaxed control  $\mu^{*n}$ solves the $n$-person game \eqref{N-person}.
  \end{enumerate}
 \end{theorem}
\begin{Proof}
  The tightness of $\big( \widebar{X}^{{n},{m},u^*},\mu^{*{n},{m}} \big)_{{m}\geq 1}$ and
  the fact that the weak limit of an extracted subsequence
  $\big( \widebar{X}^{{n},{m},u^*} \big)_{{m}\geq 1}$ solves \eqref{main_SDE_relaxed_equilibrium} with relaxed control $\mu^{*n}$
  follow from Lemma~\ref{Markov chain convergence}.
  To show \eqref{c}, 
  it suffices to prove that for all $k=0,1, \ldots ,{n}-1$ we have
\begin{equation} \label{main_approx_inf}
J\big(t_k, X^{\mu^{*n}}_{t_k}, \mu^{*n}\big) = \inf\limits_{\mu \in \mathcal{R}( [t_k,t_{k+1}] )}J\big(t_k, X^{\mu^{*n}}_{t_k}, \mu \otimes_{t_k,\Delta_n } \mu^{*n} \big).
\end{equation}
Fix any $k \in \{0,1,\ldots , {n}-1\}$,
and let $J^*_{k}$ be the infimum in the right-hand side of \eqref{main_approx_inf}.
For any $\varepsilon > 0$ there exists $\mu^{(\varepsilon )}$ such that
\[
    J^*_{k} +\varepsilon > J\big(t_k, X^{\mu^{*n}}_{t_k}, \mu^{(\varepsilon )} \otimes_{t_k,\Delta_n } \mu^{*n} \big).
\]
By Lemma~\ref{chattering lemma}, we can find
  an admissible control sequence
  $u^{{n},{m},\varepsilon} = (u_k)_{0\leq k < nm} \in \mathcal{U}^{{n},{m}}([0,T])$
such that the relaxed control representation $\mu^{{n},{m},\varepsilon}$
of $\bar{u}^{{n},{m},\varepsilon}$
converges weakly to $\mu^{(\varepsilon )}$ on $\Lambda$ as
$m$ tends to infinity. By \eqref{b},
$( \mu^{*{n},{m}} )_{m \geq 1}$ converges weakly to $\mu^{*n} $,
and therefore
$\big(\mu^{*{n},{m}}\mathbbm{1}_{[t_k,t_{k+1})^\mathsf{c}} + \mu^{{n},{m},\varepsilon}\mathbbm{1}_{[t_k,t_{k+1} )}\big)_{m \geq 1}$ converges weakly to
$\mu^{*n}\mathbbm{1}_{[t_k,t_{k+1} )^\mathsf{c}} + \mu^{(\varepsilon )} \mathbbm{1}_{[t_k,t_{k+1} )}$
on $\Lambda$ as $m \to \infty$. Then, we have
 \begin{eqnarray}
    J^*_{k} +\varepsilon  & > &J\big(t_k, X^{\mu^{*n}}_{t_k}, \mu^{(\varepsilon )} \otimes_{t_k,\Delta_n } \mu^{*n} \big)\nonumber\\
    & = &\lim\limits_{i \to \infty} J^{{n},{m_i}}\big(t_k,\widebar{X}^{{n},{m_i},u^*}_{t_k}, \bar{u}^{{n},{m_i},\varepsilon}\otimes_{t_k,\Delta_n }\bar{u}^{*{n},{m_i}}\big)
    \label{l1-1}
    \\
                          & \geq &\lim\limits_{i \to \infty} J^{{n},{m_i}}\big(t_k,\widebar{X}^{{n},{m_i},u^*}_{t_k}, \bar{u}^{*{n},{m_i}}\big)\label{main theorem approx third inequality} \\
    & = &J\big(t_k,X^{\mu^{*n}}_{t_k}, \mu^{*n}\big),
    \label{l3}
\end{eqnarray}
where $({m_i})_{{i} \geq 1}$
is an increasing sequence of integers.
\eqref{main theorem approx third inequality} is because $\bar{u}^{*{n},{m_i}}$
is solution of Problem~\eqref{problem NM-person} with Markov chain $\widebar{X}^{{n},{m_i},u^*}$,
\eqref{l1-1} and \eqref{l3}
follow from Lemma~\ref{discretized cost functional convergence},
up to extraction of a subsequence to ensure almost sure convergence.
Since $\varepsilon > 0$ is arbitrary, we conclude to \eqref{main_approx_inf}.
\end{Proof}
\medskip
\noindent
Applying Theorem~\ref{main theorem approx} requires to check the
weak convergence of a subsequence of $( \mu^{*{n},{m}} )_{{m}\geq 1}$
in $\Lambda$.
As in Corollary~\ref{main corollary}, the next result
shows that this may not be necessary
 if only the value function is concerned.
\begin{corollary}
\label{main corollary approx}
  Under Assumption~\ref{basic assumptions} and \ref{discrete_cond},
  the sequence $(J^{{n},{m}}(0,x_0,\bar{u}^{*{n},{m}}))_{{m} \geq 1}$
  admits at least one convergent subsequence.
  In addition, the limit of any such subsequence
  can be written as $J(0, x_0, \mu^{*n})$.
\end{corollary}
\begin{Proof}
  By the tightness of $( \mu^{*{n},{m}} )_{{m} \geq 1}$, we can extract a weakly convergent subsequence also denoted by $( \mu^{*{n},{m}} )_{m \geq 1}$
  whose weak limit, denoted by $\mu^{*n}$,
  is the solution to the $n$-person game \eqref{N-person}
  by Theorem~\ref{main theorem approx}, $n\geq 1$.
    By Lemma~\ref{discretized cost functional convergence} below,
    we conclude that
    $(J^{{n},{m}}(0,x_0,\bar{u}^{*{n},{m}}))_{m \geq 1}$ converges to $J(0, x_0, \mu^{*n})$.
\end{Proof}
\noindent
The next lemma has been used in the proofs of
Theorem~\ref{main theorem approx} and Corollary~\ref{main corollary approx}.

\begin{lemma}\label{discretized cost functional convergence}
  Under Assumptions~\ref{basic assumptions} and \ref{discrete_cond},
  fix $ {n} \geq 1 $ and consider a weakly convergent sequence
  $\big( \widebar{X}^{{n},{m},u},\mu^{{n},{m}}, \widebar{W}^{{n},{m},u} \big)_{{m}\geq 1}$,
  where for $m\geq 1$,
  $u^{{n},{m}} \in \mathcal{U}^{{n},{m}}([0,T])$.
  Then,  for any $k = 0,1, \ldots ,n-1$, the sequence
  $\big( J^{{n},{m}}\big(t_k,\widebar{X}^{{n},{m},u}_{t_k},\bar{u}^{{n},{m}}\big)\big)_{m\geq 1}$
  converges to $J(t_k ,X^\mu _{t_k},\mu)$ in probability as
  ${m}$ tends to infinity.
\end{lemma}
\begin{Proof}
  By the Skorokhod representation Theorem~\ref{skorohod theorem},
  Lemmas~\ref{Markov chain convergence} and \ref{equivalence stable weak},
  there is a common probability space
  $(\Omega, \mathcal{F}, \P )$ such that as $m \to \infty$, we have
\begin{equation}\label{x convergence a.s.}
\begin{cases}
\displaystyle
\lim_{{m} \to \infty}
    \sup\limits_{t \in [0,T]} \big| \widebar{X}^{{n},{m},u}_t - X^\mu _t \big|
    = 0,
\medskip
\\
\displaystyle
\lim_{{m} \to \infty}
    \left|
\int_t^T\int_U f(s,v) \mu^{{n},{m}}(ds,dv) -  \int_t^T\int_U f(s,v) \mu(ds,dv)
\right| = 0, \quad t \in [0,T],
\medskip
\\
\displaystyle
\lim_{{m} \to \infty}
\sup\limits_{t \in [0,T]} \big|\widebar{W}^{{n},{m},u}_t - W_t\big| = 0,
\end{cases}
\end{equation}
 $\P$-$a.s.$,
for $f$ any bounded random function, measurable in $t\in [0,T]$
and continuous in $u\in U$.
 Since $\big(\widebar{W}^{{n},{m},u}\big)_{t\in [0,T]}$
 is an $\big(\mathcal{F}^{{n},{m}}_t\big)_{t\in [0,T]}$-martingale
 for all $m\geq 1$,
 by Proposition~3 in \cite{briand2002robustness},
 the filtrations
 $\big(\mathcal{F}^{{n},{m}}_t\big)_{t\in [0,T]}$
 converge weakly to
 $(\mathcal{F}_t)_{t\in [0,T]}$ as $m$ tends to infinity,
 hence for all $X \in L^1(\Omega , {\cal F}, \P)$ we have the convergence
\begin{align}\label{filtration convergence in probability}
  \lim_{{m}\to \infty}
  \sup\limits_{t \in [0,T]} \mathbb{E} [X \mid \mathcal{F}^{{n},{m}}_t ]
  = \sup\limits_{t \in [0,T]}\mathbb{E} [X \mid \mathcal{F}_t ],
\end{align}
in probability.
 For any $k=0,1,\ldots , (n-1)m$, let 
\begin{align*}
Z^{n,m,u}_k :=&g\bigl(\widebar{X}^{{n},{m},u}_T, \mathbb{E}\big[\Psi(\widebar{X}^{{n},{m},u}_T) \ \big| \ \mathcal{F}^{{n},{m}}_{t^m_k}\big]\bigr) + \int_{t^m_k}^T \int_U h\bigl(s, \widebar{X}^{{n},{m},u}_s, \mathbb{E}\big[\Phi(\widebar{X}^{{n},{m},u}_s) \ \big| \ \mathcal{F}^{{n},{m}}_{t^m_k}\big], v\bigr) \mu^{{n},{m}}_s(dv) ds ,
\end{align*}
 where   $\mu^{{n},{m}}$ is
  the relaxed control representation
  of $\bar{u}^{{n},{m}}$,
  see \eqref{relaxed_control_representation}, and
$$
Z^\mu_k :=g(X^\mu _T, \mathbb{E} [\Psi(X^\mu _T) \mid \mathcal{F}_{t^m_k} ]) + \int_{t^m_k}^T \int_U h (s, X^\mu_s, \mathbb{E} [\Phi(X^\mu_s) \mid \mathcal{F}_{t^m_k} ], v ) \mu_s(dv) ds,
$$
 with
 $
 J(t^m_k,X^\mu _{t^m_k},\mu)
 =
 \mathbb{E}[Z^\mu_k \mid  \mathcal{F}_{t^m_k}]$
 and
$J^{{n},{m}}\big(t^m_k,\widebar{X}^{{n},{m},u}_{t^m_{k}},\bar{u}^{{n},{m}}\big)
 =
 \mathbb{E} \big[Z^{n,m,u}_k \mid \mathcal{F}^{{n},{m}}_{t^m_k} \big]$.
        Since convergence in 
  $L^1$ implies convergence in probability, it suffices to show that
  \[
  \lim_{{m}\to \infty}
  \mathbb{E} \big[ \big| \mathbb{E} \big[Z^{n,m,u}_k \mid \mathcal{F}^{{n},{m}}_{t^m_k} \big] - \mathbb{E}[Z^\mu_k \mid  \mathcal{F}_{t^m_k}]\big| \big] = 0.
\]
By the conditional Jensen's inequality and
Assumption~\ref{basic assumptions},
we have
\begin{align*}
  &\mathbb{E}\big[\big| \mathbb{E}\big[Z^{n,m,u}_k \mid \mathcal{F}^{{n},{m}}_{t^m_k}\big] - \mathbb{E}\big[Z^\mu_k \mid \mathcal{F}_{t^m_k}\big]
  \big| \big]  \leq \mathbb{E}\big[\big| Z^{n,m,u}_k - Z^\mu_k \big|  + | \mathbb{E}[Z^\mu_k \mid \mathcal{F}^{{n},{m}}_{t^m_k}] - \mathbb{E}[Z^\mu_k \lvert \mathcal{F}_{t^m_k}]|\big] \\
  &\leq  C \mathbb{E}\bigg[\int_{t^m_k}^T \bigl(\big|\widebar{X}^{{n},{m},u}_s - X^\mu_s\big| + | \mathbb{E}[\Phi(X^\mu_s) \mid \mathcal{F}^{{n},{m}}_{t^m_k}] - \mathbb{E}[\Phi(X^\mu_s) \mid \mathcal{F}_{t^m_k}]|\bigr) d s\\
  & \quad + \left| \int_{t^m_k}^T\bigg(\int_U h(s, X^\mu_s, \mathbb{E}[\Phi(X^\mu_s) \mid \mathcal{F}_{t^m_k}], v) \mu^{{n},{m}}_s(dv) - \int_U h(s, X^\mu_s, \mathbb{E}[\Phi(X^\mu_s) \mid \mathcal{F}_{t^m_k}], v) \mu_s(dv)\bigg) ds\right|
  \\
  &\quad +
   \big|\widebar{X}^{{n},{m},u}_T - X^\mu _T\big| + \big| \mathbb{E}\big[\Psi(X^\mu _T) \mid \mathcal{F}^{{n},{m}}_{t^m_k}\big] - \mathbb{E}\big[\Psi(X^\mu _T) \mid \mathcal{F}_{t^m_k}\big]\big|  +\big|\mathbb{E}\big[Z^\mu_k  \mid \mathcal{F}^{{n},{m}}_{t^m_k}\big] - \mathbb{E}\big[Z^\mu_k  \mid \mathcal{F}_{t^m_k}\big] \big| \bigg].
  \end{align*}
The first, third, and fourth terms in the last inequality converge
to $0$ by \eqref{x convergence a.s.},
and the fifth and sixth terms converges to $0$ by \eqref{filtration convergence in probability} and uniform boundedness. 
Similarly, by \eqref{filtration convergence in probability}
 we have
\[
\lim_{{m}\to \infty}
\mathbb{E}\big[\big| \mathbb{E}\big[\Phi(X^\mu_t) \mid \mathcal{F}^{{n},{m}}_{t^m_k}\big] - \mathbb{E} \big[\Phi(X^\mu_t) \mid \mathcal{F}_{t^m_k}\big]\big| \big]
 = 0, \qquad t\in [t^m_k,T],
\]
hence the second term tends to zero
 by the boundedness of $\Phi$ and dominated convergence.
\end{Proof}
\noindent
\noindent
{\em Remark.}
 In addition to the dependence of
  $h(s,X^{\mu}_s,\mathbb{E}_t [\Phi(X^{\mu}_s)], u)$
 and
 $g(X^{\mu}_T, \mathbb{E}_t[\Psi(X^{\mu}_T])$
 on the mean-field term, time inconsistency of a control problem can also be caused by the dependence of $h$ and $g$ on initial time and initial state
 $t$ and $X^{\mu}_t$, i.e.
 \[
    J(t,X^{\mu}_t,\mu) = \mathbb{E}_t\left[
     g\big(t,X^{\mu}_t,X^{\mu}_T, \mathbb{E}_t[\Psi(X^{\mu}_T)]\big) +
    \int_t^T \int_U h\big(t,s,X^{\mu}_t,X^{\mu}_s, \mathbb{E}_t [\Phi(X^{\mu}_s)], v\big)\mu_s(dv) ds\right],
\]
which admits the discretization
\begin{eqnarray}
  \label{cost_functional_discrete_extended}
  J^{{n},{m}}\bigl(t,\widebar{X}^{{n},{m},u}_{t},\bar{u}^{{n},{m}}\bigr) & = &\mathbb{E}_t^{{n},{m}}\bigg[
    g\bigl(t,\widebar{X}^{{n},{m},u}_{t},\widebar{X}^{{n},{m},u}_{T}, \mathbb{E}_t^{{n},{m}}\big[\Psi\big(\widebar{X}^{{n},{m},u}_{T}\big)\big]\bigr)
    \\
    \nonumber
    & & \quad +
\int_t^T h\bigl(t,s,\widebar{X}^{{n},{m},u}_{t},\widebar{X}^{{n},{m},u}_s, \mathbb{E}_t^{{n},{m}}\big[\Phi\big(\widebar{X}^{{n},{m},u}_s\big) \big], \bar{u}^{{n},{m}}_s\bigr) ds
    \bigg],
\end{eqnarray}
where $\E^{{n},{m}}_t[ \ \! \cdot \ \! ] = \E\big[ \ \! \cdot \! \mid \! \mathcal{F}_t^{{n},{m}}\big]$.
We note that under additional uniform continuity and Lipschitz continuity assumptions on
$h(t,s,\xi ,x,y,u)$,
and
$g(t,\xi ,x,y)$ in initial time $t$ and initial state $\xi $ respectively,
the analysis
of Theorem~\ref{main theorem},
Corollary~\ref{existence}
and Theorem~\ref{main theorem approx}
can be extended to the setting of \eqref{cost_functional_discrete_extended},
by replacing \eqref{hamiltonian} with the Hamiltonian
\begin{equation*}
H(t,s,\xi ,x,y,\mu,p) = p\int_U b (s,x,v ) \mu(dv) - \int_U h (t,s,\xi ,x,y, v ) \mu(dv).
\end{equation*}
The proofs of Section~\ref{sec:characterization_nperson} remain unchanged because the spike perturbation
$\mu \otimes_{t,\varepsilon} \mu^{*n}$ does not affect the initial state $X^{\mu^{*n} }_{t^m_k}$.
The main changes to Section~\ref{sec:convergence_nperson} are in Theorem~\ref{main theorem},
where the bound \eqref{main theorem last inequality} on
\begin{equation*}
    \mathbb{E}\left[\int_0^T\big| H\bigl(\floor{t}_{n},t,X^{{\mu}^n }_{\floor{t}_{n}},X^{{\mu}^n }_t,\mathbb{E}_{\floor{t}_{n}}\big[ \Phi\big( X^{{\mu}^n }_t \big) \big],\mu^{*n} _t,p^{\mu^{*n} }_{\floor{t}_{n},t}\bigr) - H\bigl(t,t,X_t,X_t,\Phi( X_t ),\mu^*_t,p^{\mu^*}_{t,t}\bigr)\big| dt\right]
\end{equation*}
now contains two additional terms
\begin{equation*}
\int_0^T \mathbb{E} \biggl[
    \lvert X_{\floor{t}_{n}} - X_t \rvert +
    \int_U \bigl| h\bigl(\floor{t}_{n},t,X_{t},X_t,\mathbb{E}_{t}[\Phi(X_t)], v\bigr) - h\bigl(t,t,X_t,X_t,\mathbb{E}_t[\Phi(X_t)], v\bigr)\bigr| \mu^{*n} _t(dv)  \biggr]dt,
\end{equation*}
which converge to $0$ by noting the uniform continuity of $h$ on initial time and the continuity property of SDE.
The proofs in Section~\ref{sec:markovchain_approx}, particularly Lemma~\ref{discretized cost functional convergence}, can be modified similarly.

\subsection{Numerical results}
\label{sec:numerical}
In this section we present numerical illustrations based on
Theorem~\ref{main theorem approx}.
Assume that $K$ is the bounding constant in Assumption~\ref{basic assumptions},
and let $p^{{n},{m}}(y; t^m_k, x, u)$ denote the transition probability
of $(X^{{n},{m},u}_{k+1})_{0\leq k < nm}$,
$x \in \sqrt{\Delta}_{{n},{m}}\mathbb{Z}$,
$u \in \mathcal{U}^{{n},{m}}([0,T])$.
   As in \S~4 of \cite{fischer2007discretisation},
   Assumption~\ref{discrete_cond}
   is satisfied using a trinomial tree constructed as
\begin{equation*}
    p^{{n},{m}} (y; t^m_k, x, u )=
\begin{cases}
  \displaystyle
  \frac{\sqrt{\Delta}_{{n},{m}}}{2K} b(t^m_k, x, u_k) +\frac{1}{2K^2} \sigma^2(t^m_k, x), &
 y = x + K\sqrt{\Delta}_{{n},{m}} ,
\medskip
\\
\displaystyle
  -\frac{\sqrt{\Delta}_{{n},{m}}}{2K} b(t^m_k, x, u_k) + \frac{1}{2K^2} \sigma^2(t^m_k, x), &
 y = x - K\sqrt{\Delta}_{{n},{m}},
\medskip
\\
\displaystyle
  1 - \frac{1}{K^2} \sigma^2(t^m_k, x), &
 y = x ,
\medskip
\\
0 ,                     & \text{ otherwise.}
\end{cases}
\end{equation*}
\noindent
 We consider the following numerical implementation of Theorem~\ref{main theorem approx}.
\begin{enumerate}[(i)]
      \item For each time $t^m_k$, initialize the nodes ${\cal Y}_k := \big\{ H_{{n},{m}}(x_0) + j K \sqrt{\Delta}_{{n},{m}} \ : \ - k \leq j \leq k \big\}$.
      \item Starting from $t^m_{({n}-1){m}}$,
        solve Problem~\eqref{problem NM-person}
        for every initial value $x \in {\cal Y}_{({n}-1){m}}$
        at time $t^m_{({n}-1){m}}$.
      \item Repeat (ii) recursively at times $t^m_{({n}-2){m}}, \dots, t^m_{m}, t^m_0$.
\end{enumerate}
However, solving Problem~\eqref{problem NM-person} can still be computationally expensive for large ${m}$ because we need to optimize $1+3+3^2+\dots+3^{{m}-1} = (3^{m}-1)/2$
controls at each node $x \in {\cal Y}_{t_k }$,
 $k   = {n}-1, {n}-2, \ldots, 0$.
If the function \eqref{cost_functional_discrete_extended} does not depend on a
mean-field term then for each node $x \in {\cal Y}_{t_k }$,
$k = {n}-1, {n}-2, \ldots, 0$, the optimization problem
\begin{equation}\label{sub_problem}
    \inf\limits_{u \in \mathcal{U}^{{n},{m}}( [0,T ] )} J^{{n},{m}} (t_k,x , \bar{u} \otimes_{t_k,t_{k+1} } \bar{u}^{{n},{m}})
\end{equation}
 can be solved using dynamic programming,
which reduces the number of parameters
to be optimized
from exponential $(3^{m}-1)/2$ to polynomial $1 + 3 + 5 + \cdots + ( 2({m}-1) + 1 ) = {m}^2$ at every node $x \in {\cal Y}_{t_k }$, $k = {n}-1, {n}-2, \ldots, 0$.

\medskip

To solve \eqref{sub_problem} using
dynamic programming at each time $t^m_l$
with $t_k \le t^m_l < t_{k+1}$,
we need to access the optimal control on $[t^m_{l+1} , T]$ and
calculate $J^{{n},{m}}(t^m_l,x , u)$,
which involves a calculation from time $t^m_l$ to time $T$.
The complexity of the algorithm
can be reduced in case \eqref{cost_functional_discrete_extended}
takes the particular form
\begin{eqnarray}
  \label{specific_func}
  \lefteqn{
    \! \! \! \! \! \! \! \! \! \! \!
     J^{{n},{m}}\big(t,\widebar{X}^{{n},{m},u}_{t},\bar{u}^{{n},{m}}\big)
 = g\bigl(t,\widebar{X}^{{n},{m},u}_{t},\mathbb{E}^{{n},{m}}_t [g_1(\widebar{X}^{{n},{m},u}_{T})],
    \ldots, \mathbb{E}^{{n},{m}}_t [g_p (\widebar{X}^{{n},{m},u}_{T})]\bigr)
  }
  \\
  \nonumber
     & & + \int_t^T h\bigl(t,s,\widebar{X}^{{n},{m},u}_{t}, \mathbb{E}^{{n},{m}}_t [h_1(\widebar{X}^{{n},{m},u}_{s}, u_s)],
  \ldots, \mathbb{E}^{{n},{m}}_t [h_q(\widebar{X}^{{n},{m},u}_{s}, u_s)]\bigr) ds,
  \nonumber
\end{eqnarray}
 from which we have 
\begin{align}
&J^{{n},{m}}(t^m_l,x , \bar{u}^{{n},{m}})  =g\bigl(t^m_l,x,\mathbb{E}^{{n},{m}}_{t^m_l} \big[ \mathbb{E}^{{n},{m}}_{t^m_{l+1}} \big[g_1\big(\widebar{X}^{{n},{m},u}_{T}\big)\big] \big],
\ldots, \mathbb{E}^{{n},{m}}_{t^m_l} \big[ \mathbb{E}^{{n},{m}}_{t^m_{l+1} } \big[g_p \big(\widebar{X}^{{n},{m},u}_{T}\big)\big] \big]\bigr)
    \nonumber
    \\
& \quad+\int_{t^m_l}^{t^m_{l+1}} h\bigl(t^m_l,s,x, \mathbb{E}^{{n},{m}}_{t^m_l} \big[h_1\big(\widebar{X}^{{n},{m},u}_s, \bar{u}^{{n},{m}}_s\big)\big],
    \ldots, \mathbb{E}^{{n},{m}}_{t^m_l} \big[h_q\big(\widebar{X}^{{n},{m},u}_s, \bar{u}^{{n},{m}}_s\big)\big]\bigr) ds
    \nonumber
    \\
&  \quad+ \int_{t^m_{l+1} }^{T} h\bigl(t^m_l,s,x, \mathbb{E}^{{n},{m}}_{t^m_l} \big[ \mathbb{E}^{{n},{m}}_{t^m_{l+1}} \big[h_1\big(\widebar{X}^{{n},{m},u}_s, \bar{u}^{{n},{m}}_s\big)\big] \big],
\ldots, \mathbb{E}^{{n},{m}}_{t^m_l} \big[ \mathbb{E}^{{n},{m}}_{t^m_{l+1} } \big[h_q\big(\widebar{X}^{{n},{m},u}_s, \bar{u}^{{n},{m}}_s\big)\big] \big]\bigr) ds.
    \nonumber
    \\
    \label{specific_func2}
\end{align}
 In this case it suffices to maintain an array for
the values of $h_1, \ldots, h_q, g_1, \ldots, g_p$
at $\widebar{X}^{{n},{m},u}_{s}, \bar{u}^{{n},{m}}_s$  in order to solve \eqref{specific_func2}
at time $t^m_l$,
which involves calculations from time $t^m_l$ to time $t^m_{l+1}$,
instead of from $t^m_l$ to $T$.
This method is applied to the quadratic and quartic cost functions
examples $\E^{n,m}_t\big[\big(\widebar{X}^{{n},{m},u}_{T} - \widebar{X}^{{n},{m},u}_{t}\big)^2\big]$ and
$\E^{n,m}_t \big[\big(\widebar{X}^{{n},{m},u}_{T} - \widebar{X}^{{n},{m},u}_{t}\big)^4\big]$
below,
 however not all cost functions satisfy \eqref{specific_func},
 e.g. $\E^{n,m}_t[1 / (\widebar{X}^{{n},{m},u}_{T} + \widebar{X}^{{n},{m},u}_{t} )] $
 cannot be written in that form.

\subsubsection{Linear-quadratic control problem}
\begin{figure}[H]
  \centering
 \begin{subfigure}[b]{0.45\textwidth}
    \includegraphics[width=1\linewidth]{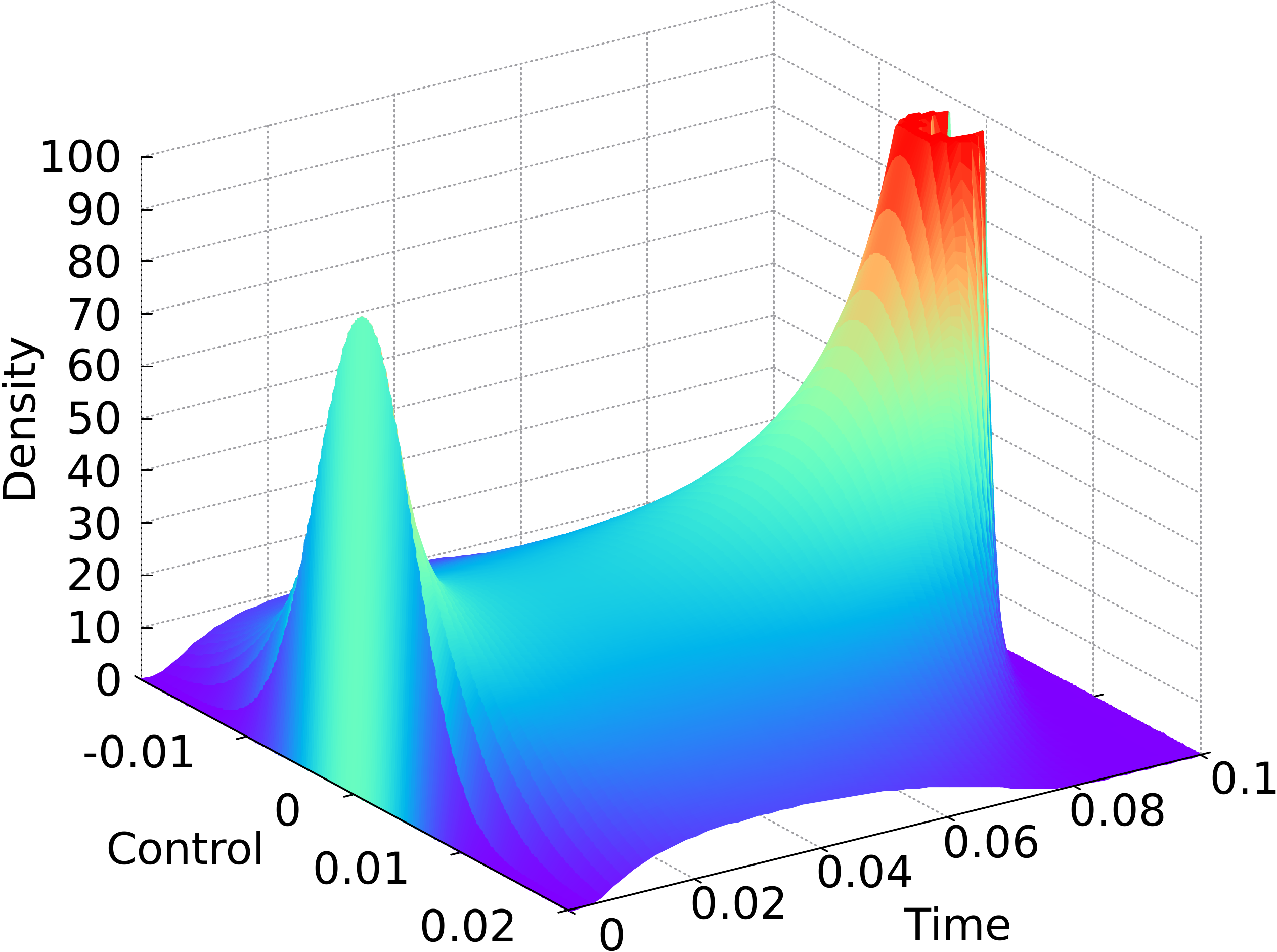}
    \caption{Equilibrium controls.}
    \label{fig:density_quadratic_equilibrium_20person-a}
 \end{subfigure}
  \begin{subfigure}[b]{0.45\textwidth}
    \includegraphics[width=1\linewidth]{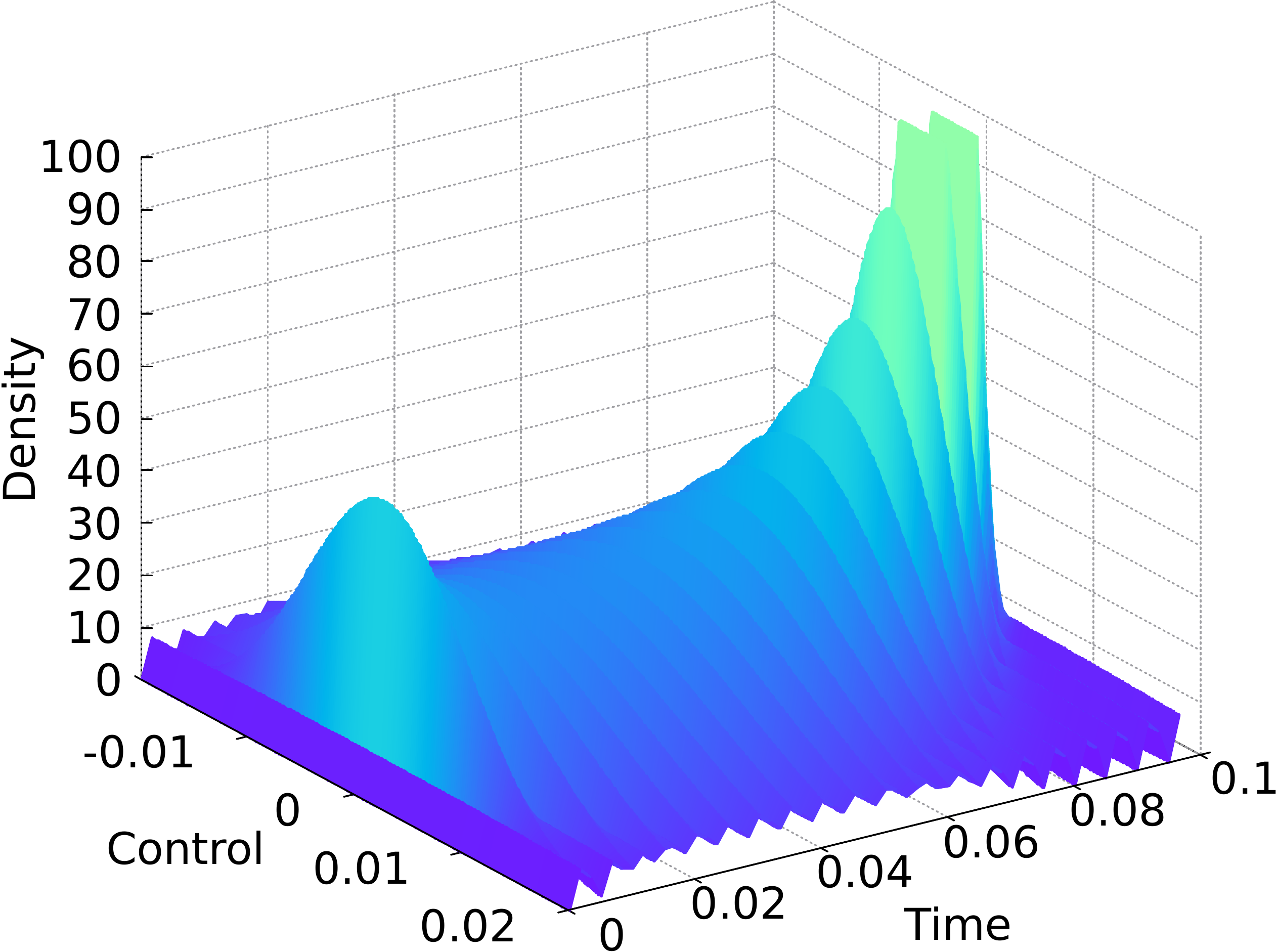}
    \caption{Solution of $20$-person game.}
    \label{fig:density_quadratic_equilibrium_20person-b}
 \end{subfigure}
    \caption{Comparison between the equilibrium control and the 20-person game solution.}
\label{fig:density_quadratic_equilibrium_20person}
\end{figure}

We first check the numerical application of
Theorems~\ref{main theorem} and \ref{main theorem approx}
to a linear-quadratic control problem which admits an analytic solution,
 see \cite{bjork2010general} and \cite{djehiche2016characterization},
 allowing us to evaluate the performance of our numerical scheme.
Here, the state of the system is driven by the SDE
\begin{equation}
  \label{LQ_SDE}
  \left\{
  \begin{array}{l}
    \displaystyle dX^{\xi,\mu}_{t,s} = \left(
    aX^{\xi,\mu}_{t,s} + c \int_U v \mu_s(dv)\right) ds + \sigma dW_s,
     \quad t \leq s \leq T,
\medskip
\\
X^{\xi,\mu}_{t,t} = \xi,
  \end{array}
  \right.
\end{equation}
where $a, c, \sigma \in \mathbb{R}$,
 with the cost functional
 \begin{equation}
   \label{LQ_cost_functional}
  J(t,\xi ,\mu) =
  \frac{\gamma}{2}
  \mathbb{E}\big[ \big(X^{\xi ,\mu}_{t,T} - \xi \big)^2 \big]
    + \frac{1}{2}\mathbb{E}\left[\int_t^T \int_U v^2 \mu_s(dv) ds\right],
\end{equation}
 where $\gamma > 0$,
 $g(t,\xi ,x,y) = \gamma (x-\xi)^2/2$
 and
 $h(t,s,\xi ,x,y,u) = u^2/2$,
 in the framework of \eqref{cost_functional_discrete_extended}.
Extending the solution technique of \cite{djehiche2016characterization}
from the strict control space to the relaxed control space
by replacing Theorem~1 therein with Theorem~\ref{MP TC} above,
it can be shown that \eqref{LQ_cost_functional} admits
 a strict equilibrium control represented as
\begin{equation}
  \label{mun1}
 \mu^*_t (dv) = \delta_{u^*_t}(dv)
  := \delta_{c(\beta (t)-\alpha (t)) X^{\mu^*}_t}(dv),
\end{equation}
 where
$$
 X^{\mu^*}_t :=
X^{x_0 , \mu^*}_{0,t}=
 x_0 \Gamma(0,t) + \sigma \int_0^t \Gamma(s, t) dW_s,
\quad
\Gamma(t, s) = \exp\left(\int_t^s \big(
a + c^2 (\beta (r) - \alpha (r) ) \big)
dr\right),
$$
 and the functions $\alpha (t)$, $\beta (t)$ are defined by
$$
\beta (t) =\gamma e^{a(T-t)}, \quad
\alpha (t) = \gamma \frac{
 \exp\big(2a(T-t) + c^2\int_t^T \beta (s) ds\big)}{ 1 + \gamma  c^2 \int_t^T \exp\big(
2a(T-s) + c^2 \int_s^T \beta (r) dr\big) ds}.
$$
\begin{prop}
  Let $k\in \{1,\ldots ,n\}$.
  The solution of the $n$-person game
  \begin{equation}
\nonumber 
J\bigl(t_{ k-1 },X^{\mu^{*n} }_{t_{ k-1 }}, \mu^{*n} \bigr) = \inf\limits_{\mu \in \mathcal{R}( [t_{ k-1 },t_k] )} J\bigl(t_{ k-1 },X^{\mu^{*n} }_{t_{ k-1 }}, \mu \otimes_{t_{ k-1 },\Delta_n} \mu^{*n} \bigr)
  \end{equation}
  is given by the strict equilibrium control represented as
\begin{equation}
  \label{mun0}
  \mu^{*n} _t (dv)
  = \delta_{u^{*n}_t} (dv)
  := \delta_{c  \beta_{n} (t) X^{\mu^{*n} }_{t_{ k-1 }} - c  \alpha_{n} (t) X^{\mu^{*n} }_{t}} (dv),
  \qquad t_{ k-1 } < t \leq t_k,
\end{equation}
where
\begin{equation}
   \label{fjklds}
\left\{
  \begin{array}{ll}
    \displaystyle
    \alpha_{n} (t)&=     \displaystyle\frac{2a\alpha_{n} (t_k) e^{-2a (t-t_k)}}{
        2a + c^2\alpha_{n} (t_k) ( e^{-2a (t-t_k)} -1)
    },
    \medskip
 \\
    \displaystyle
    \beta_{n} (t) &=     \displaystyle
    \beta_{n} (t_k)\exp\left(a(t_k-t) -  c^2 \int_t^{t_k} \alpha_{n} (s) ds\right),
    \qquad t_{ k-1 } < t \leq t_k,
  \end{array}
  \right.
\end{equation}
with the terminal conditions
\begin{equation}
  \label{term}
  \alpha_n (t_k) = \gamma e^{a(T-t_k)}\prod\limits_{l=k}^{{n}-1}\bigg( \Gamma_n (t_l, t_{l+1}) + c \int_{t_l}^{t_{l+1}} \Gamma_n (s, t_{l+1})\beta_n(s) ds \bigg),
  \ \ \beta_{n} (t_k) = \gamma e^{a(T-t_k)},
  \end{equation}
where
 $\Gamma_n (t, s) = \exp\left(\int_t^s ( a - c^2 \alpha_n(r) ) dr\right)$.
\end{prop}
\begin{Proof}
      We work by backward induction, starting from $k=n$.
  In addition to proving \eqref{mun0}-\eqref{term}, we also show that
    \begin{equation}
        \label{E XT}
        \E_{t_k} \big[ X^{\mu^{*n} }_T \big] = X^{\mu^{*n} }_{t_k} \prod\limits_{l=k}^{n-1}\left( \Gamma_n (t_l, t_{l+1}) + c \int_{t_l}^{t_{l+1}} \Gamma_n (s, t_{l+1})\beta_n(s) ds \right).
    \end{equation}
 The corresponding adjoint equation \eqref{TC adjoint discrete} can be written as
\begin{equation}
  \label{LQ bsde}
  \left\{
  \begin{array}{l}
    \displaystyle
    dp^{k-1,\mu^{*n} }_{t_{k-1},t} = -ap^{k-1,\mu^{*n} }_{t_{k-1},t} dt + q^{k-1,\mu^{*n} }_{t_{k-1},t} dW_t, \qquad t_{k-1} \leq t \leq T ,
\medskip
 \\
    \displaystyle
    p^{k-1,\mu^{*n} }_{t_{k-1},T} = \gamma \big(X^{\mu^{*n} }_{t_{k-1}} - X^{\mu^{*n} }_T\big),
  \end{array}
  \right.
\end{equation}
 with solution given by
\begin{equation}
\nonumber 
    p^{k-1,\mu^{*n} }_{t_{k-1},t} = \gamma e^{a(T-t)}\big(X^{\mu^{*n} }_{t_{k-1}} - \E_t\big[ X^{\mu^{*n} }_T\big] \big),
\end{equation}
hence by \eqref{E XT} we have
\begin{equation}
    \label{lq bsde terminal cond tk}
    p^{k-1,\mu^{*n} }_{t_{k-1},t_k} = \gamma e^{a(T-t_k)}X^{\mu^{*n} }_{t_{k-1}} - \gamma e^{a(T-t_k)}
    X^{\mu^{*n} }_{t_k}
    \prod\limits_{l=k}^{n-1}\left( \Gamma_n (t_l, t_{l+1}) + c \int_{t_l}^{t_{l+1}} \Gamma_n (s, t_{l+1})\beta_n(s) ds \right).
\end{equation}
 Next, we look for the solution of the form
\begin{equation}
  \label{dflksf1}
  p^{k-1,\mu^{*n} }_{t_{k-1},t} = \beta_n (t) X^{\mu^{*n} }_{t_{k-1}} - \alpha_n (t) X^{\mu^{*n} }_{t},
 \qquad t_{k-1} < t \leq t_k .
  \end{equation}
 By It\^o's lemma and \eqref{LQ_SDE}, we have
\begin{eqnarray*}
  dp^{k-1,\mu^{*n} }_{t_{k-1},t} & = &
  \beta_n' (t) X^{\mu^{*n} }_{t_{k-1}} dt
 - \alpha_n' (t) X^{\mu^{*n} }_{t} dt
 - \alpha_n (t) dX^{\mu^{*n} }_{t}
 \\
 & = &
  \beta_n' (t) X^{\mu^{*n} }_{t_{k-1}} dt
 - \alpha_n' (t) X^{\mu^{*n} }_{t} dt
 - \alpha_n (t)
 \left( aX^{\mu^{*n}}_t + c \int_U v \mu^{*n}_t(dv)\right) dt
  -  \sigma \alpha_n (t)
dW_t,
 \end{eqnarray*}
 and comparing the resulting coefficients in `$dt$' and `$dW_t$' with \eqref{LQ bsde}, we obtain
\begin{equation}\label{ode terms}
  (\alpha_n' (t) + 2a\alpha_n (t) )X^{\mu^{*n} }_{t}
  + c \alpha_n (t) \int_U v \mu^{*n}_t(dv) dt
 = (\beta_n' (t) + a\beta_n (t) )X^{\mu^{*n} }_{t_{k-1}},
 \quad t_{k-1} < t \leq t_k,
\end{equation}
and $q^{k-1,\mu^{*n} }_{t_{k-1},t} = -\sigma\alpha_n (t)$.
 By \eqref{hamiltonian}, the Hamiltonian of this system is
\[
H(t, x, \mu, p) := axp + c p \int_U v \mu(dv) - \frac{1}{2} \int_U v^2 \mu(dv).
\]
Due to the concavity of $H(t, x, \mu, p)$,
the optimality necessary condition \eqref{MP k H compare N person} in
Theorem~\ref{MP N person} becomes sufficient,
see i.e. Theorem~3.5.2 in \cite{yong1999stochastic}
and Theorem~4.1 in \cite{andersson2011},
and it yields
$\mu^{*n} _t (dv) = \delta_{c p^{k-1,\mu^{*n} }_{t_{k-1},t}} (dv)$
 on $(t_{k-1}, t_k]$ after maximizing
      $H(t, X^{\mu^{*n}}_t, \cdot , p^{k-1,\mu^{*n} }_{t_{k-1},t} )$,
      which shows \eqref{mun0}.
      Next, plugging \eqref{mun0} into \eqref{ode terms}
      and identifying the coefficients in `$X^{\mu^{*n} }_{t_{k-1}}$' and `$X^{\mu^{*n} }_{t}$', we obtain
\[
\left\{
  \begin{array}{l}
    \displaystyle
    \alpha_n' (t) + 2a\alpha_n (t) - c^2(\alpha_n (t))^2 = 0,
\medskip
 \\
    \displaystyle
    \beta_n' (t) + (a - c^2\alpha_n (t) )\beta_n (t) = 0,
    \qquad
    t_{k-1} < t \leq t_k,
  \end{array}
  \right.
\]
which yields \eqref{fjklds}, while the terminal conditions \eqref{term}
are obtained by a comparison of \eqref{lq bsde terminal cond tk} and \eqref{dflksf1}.
Regarding \eqref{E XT}, we have
\begin{align*}
  \E_{t_{k-1}} \big[ X^{\mu^{*n} }_T \big] & =
  \E_{t_{k-1}} \big[ \E_{t_k} \big[ X^{\mu^{*n} }_T \big] \big]
  \\
  & =
  \E_{t_{k-1}} \left[X^{\mu^{*n} }_{t_k} \prod\limits_{l=k}^{n-1}\left( \Gamma_n (t_l, t_{l+1}) + c \int_{t_l}^{t_{l+1}} \Gamma_n (s, t_{l+1})\beta_n(s) ds \right)\right]
   \\
  & =
   \E_{t_{k-1}} \big[X^{\mu^{*n} }_{t_k} \big]
   \prod\limits_{l=k}^{n-1}\left( \Gamma_n (t_l, t_{l+1}) + c \int_{t_l}^{t_{l+1}} \Gamma_n (s, t_{l+1})\beta_n(s) ds \right)
   \\
                                                & = X^{\mu^{*n} }_{t_{k-1}} \prod\limits_{l=k-1}^{n-1}\left( \Gamma_n (t_l, t_{l+1}) + c \int_{t_l}^{t_{l+1}} \Gamma_n (s, t_{l+1})\beta_n(s) ds \right),
\end{align*}
where the last equality is obtained by solving the linear SDE \eqref{LQ_SDE}
using \eqref{mun0}.
 Finally, assuming that \eqref{mun0}-\eqref{E XT} hold at the rank $k$,
 we repeat the above argument to show that they hold at the rank $k-1$.
\end{Proof}
\noindent
In Figure~\ref{fig:density_quadratic_equilibrium_20person} we compare
the actual probability density of the equilibrium control $u^*_t$ given by
\eqref{mun1} to the $20$-person game solution $u^{*20}_t$
obtained from \eqref{mun0} for $t\in [0,T]$ with $T=0.1$.

\medskip

\begin{figure}[H]
  \centering
  \includegraphics[width=\textwidth]{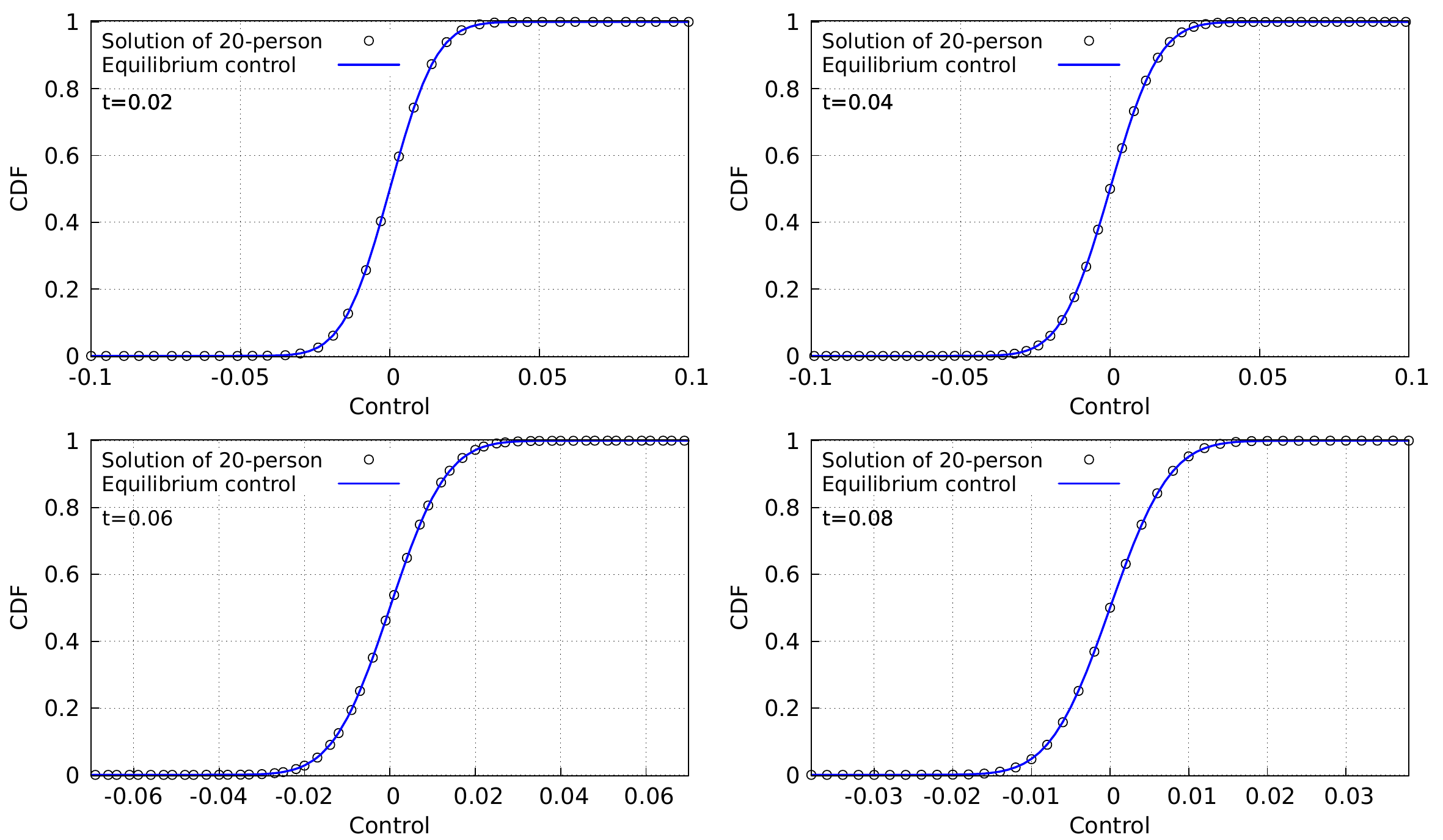} 
  \caption{CDF comparison between $\mu^*$ and the $20$-person game solution.}
  \label{fig:quadratic_equilibrium_20person}
\end{figure}

 \noindent
  In Figure~\ref{fig:quadratic_equilibrium_20person} 
  we check the convergence in distribution of $u^{*n}_t$ in
\eqref{mun0} to $u^*_t$ in \eqref{mun1} by comparing the CDFs
of $u^*_t$ and $u^{*n}_t$ with $n=20$ at times $t=0.02, 0.04, 0.06, 0.08$,
with $a=c =\sigma=\gamma=1$ and $x_0=0$. 

\subsubsection*{Numerical approximation of the $n$-person game solution}
\begin{figure}[H]
  \centering
  \includegraphics[width=\textwidth]{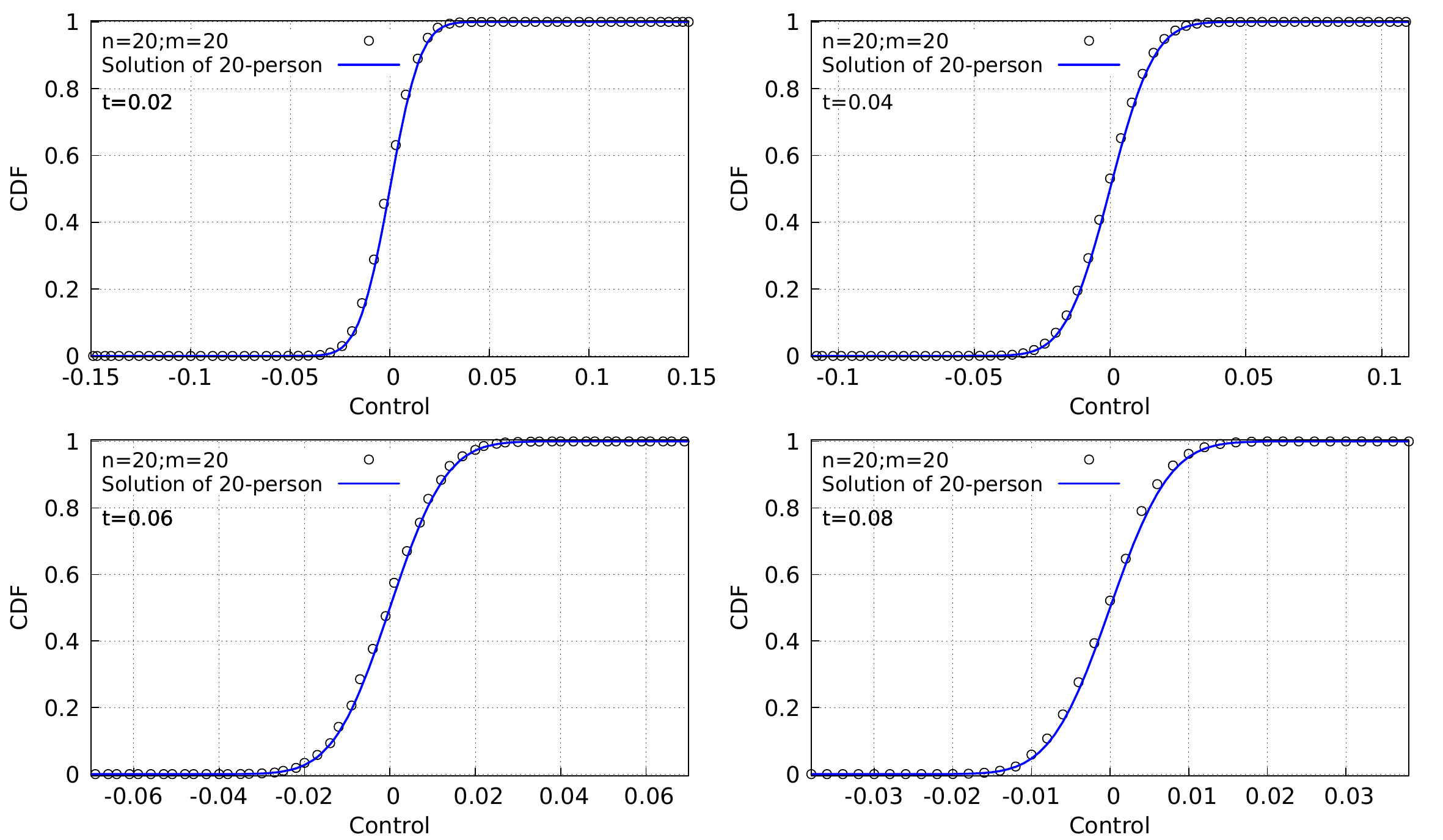}
  \caption{CDF Comparison between
    the $20$-person game solution 
    and the numerical solution.} 
\label{fig:quadratic_cdf_diff_time}
\end{figure}

\noindent
To assess the weak convergence of controls stated in Theorem~\ref{main theorem} and \ref{main theorem approx}, in Figure~\ref{fig:quadratic_cdf_diff_time}
we compare the closed form CDFs of
$u^{*n}_t$ obtained from \eqref{mun0}
to the numerical solution $\bar{u}^{*{n},{m}}_t$ of Problem~\eqref{problem NM-person}
with $n = m = 20$ at times $t =0.02, 0.04, 0.06, 0.08$,
and $U=[-10,10]$, by truncating
$b(t,x,u)$, $h(t,x,y,u)$, $g(x,y)$ up to $K$.

\begin{figure}[H]
\centering
   \includegraphics[width=0.65\linewidth]{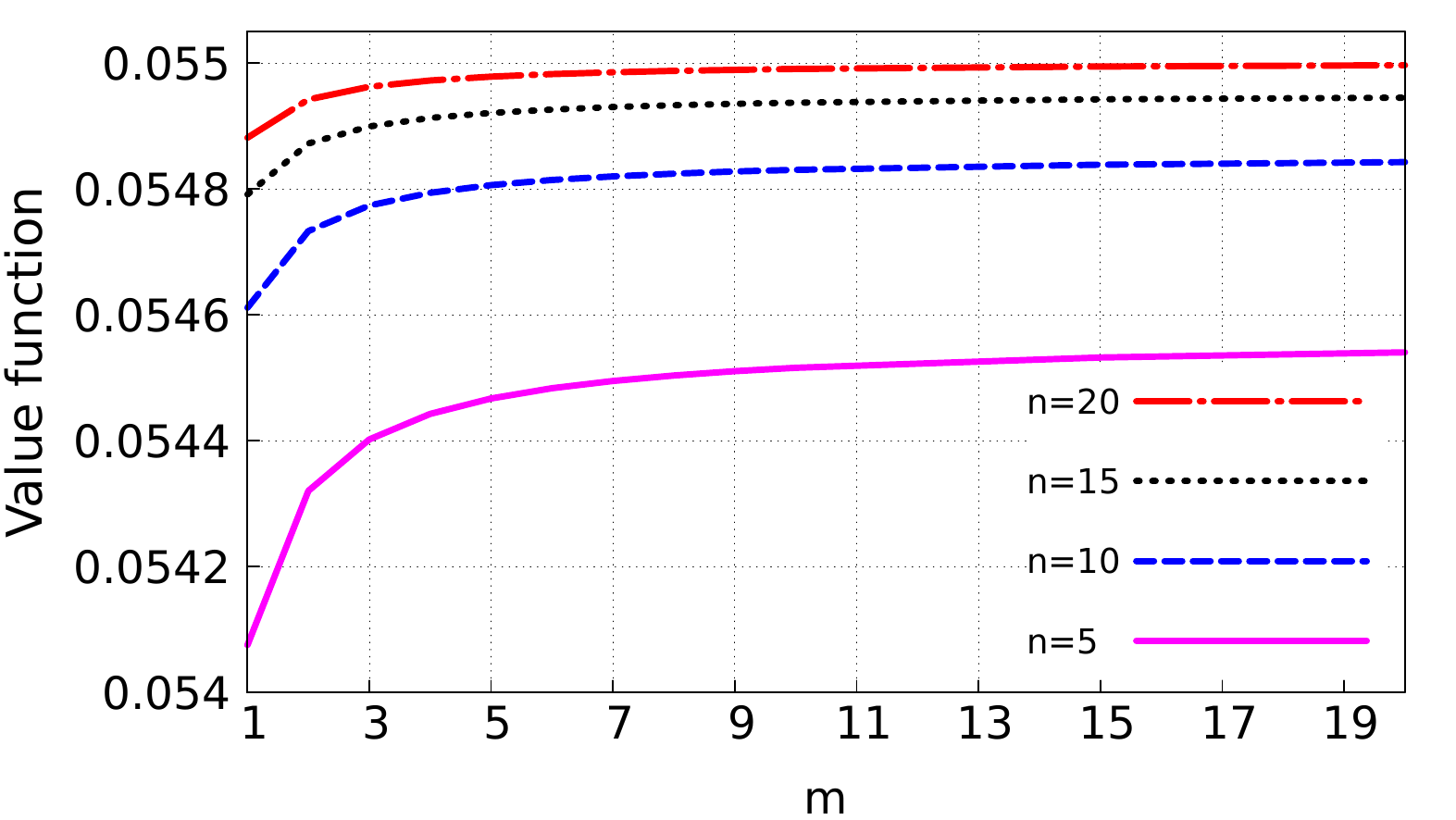} 
   \caption{Comparison of value functions.}
   \label{fig:quadratic_value_func(a)}
\end{figure}

\noindent
 In Figure~\ref{fig:quadratic_value_func(a)},
 we compare the value functions $J^{n,m}(0, x_0, \bar{u}^{*n,m})$
 with
${n} = 5, 10, 15, 20$ and ${m} \in \{1, \ldots , 20\}$.

\begin{figure}[H]
  \centering
  \includegraphics[width=0.65\linewidth]{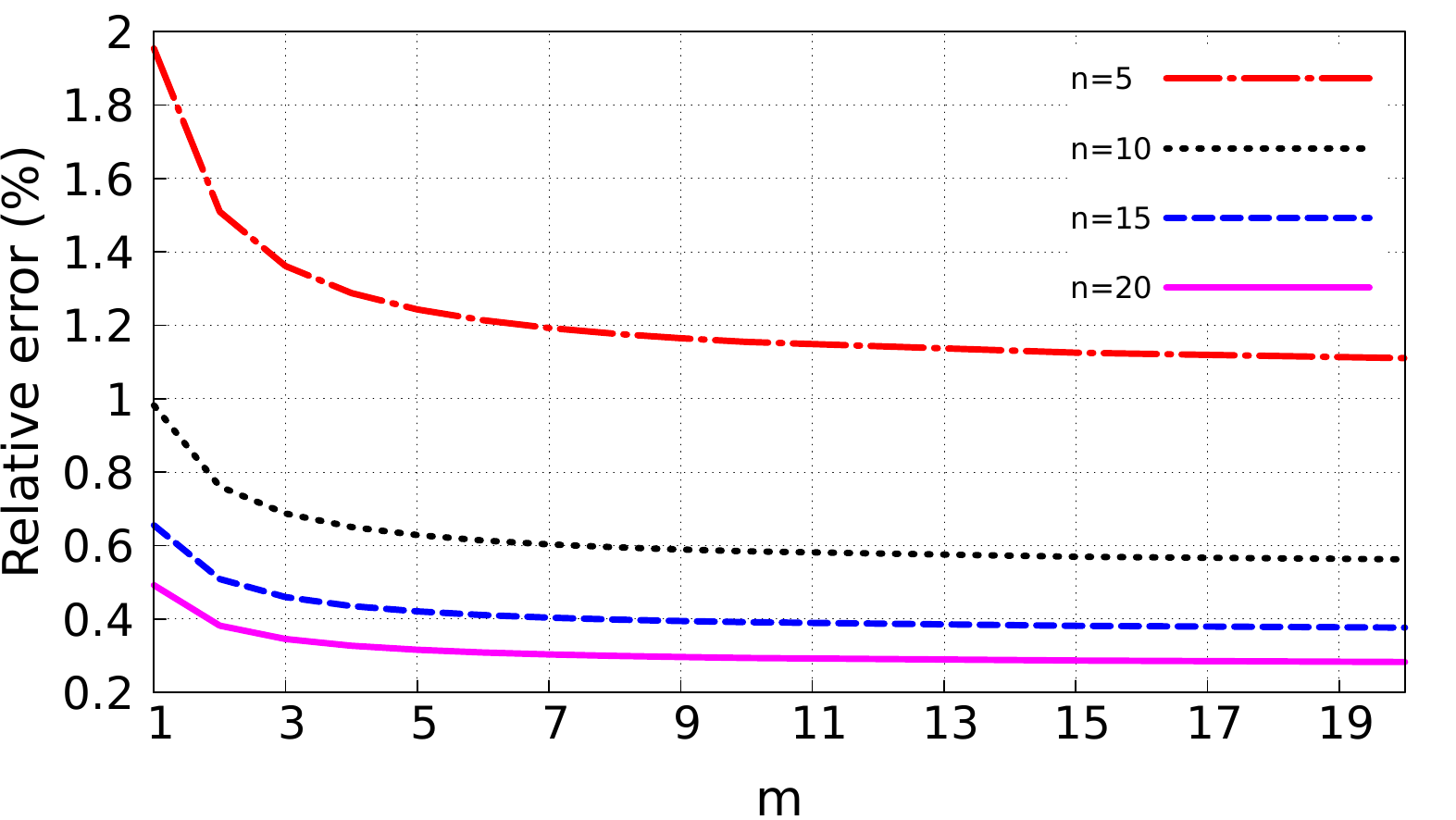}
   \caption{Comparison of relative errors in percentage.}
   \label{fig:quadratic_value_func(b)}
\end{figure}

\noindent
In Figure~\ref{fig:quadratic_value_func(b)}, we compare the
relative errors of the value function $J^{n,m}(0, x_0, \bar{u}^{*n,m})$
with respect to $J(0, x_0, \mu^*)$.

\subsubsection{Linear-quartic control problem}
\begin{figure}[H]
\centering
   \includegraphics[width=0.65\linewidth]{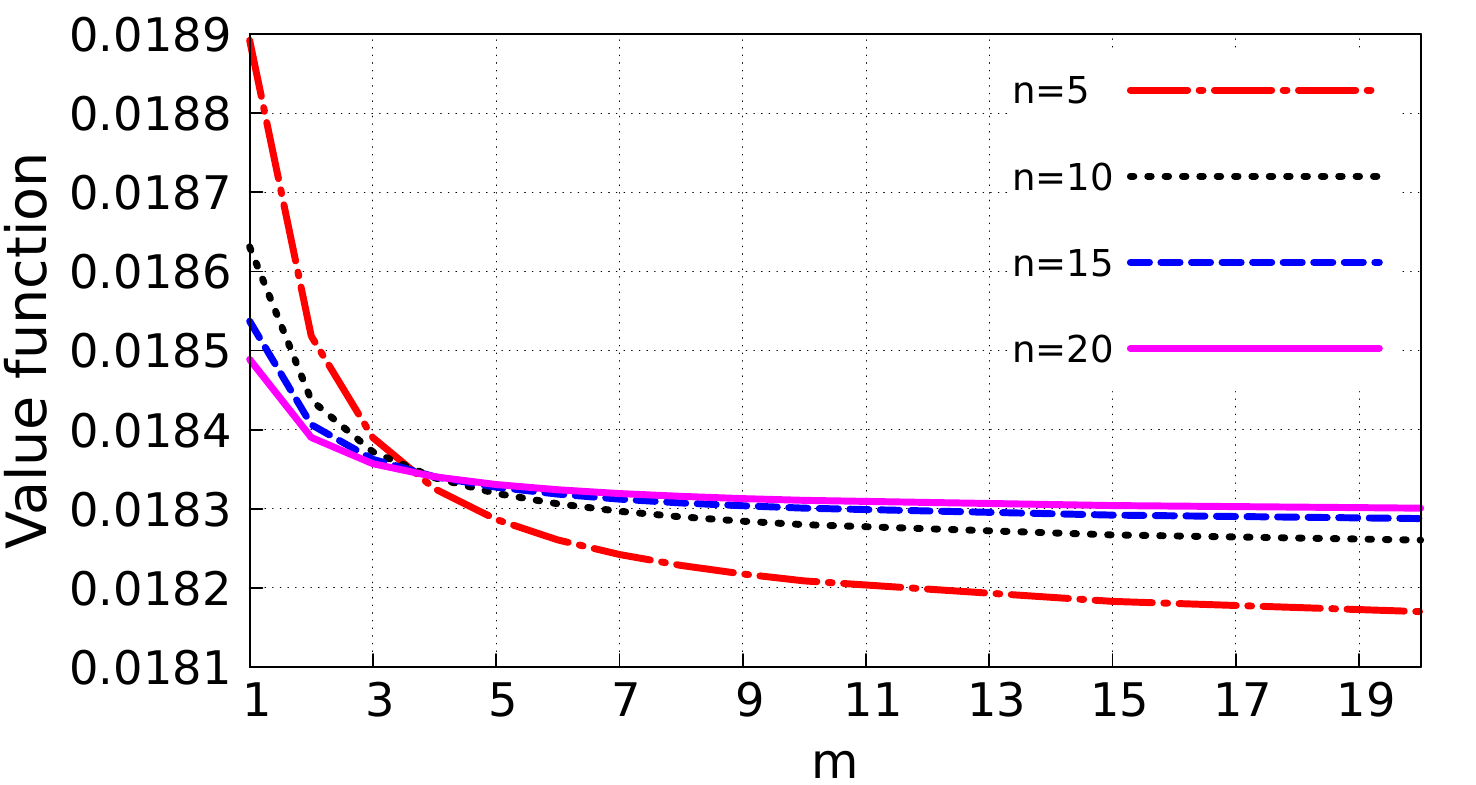}
   \caption{Comparison of value functions.}
   \label{fig:quartic_value_func(a)}
\end{figure}

Here, we apply our solution algorithm to the problem
\begin{equation}
  \nonumber 
  J(t,\xi,\mu) =
    \frac{\gamma}{2} \mathbb{E} [(X^{\xi,\mu}_{t,T} - \xi)^4]
    + \frac{1}{2} \mathbb{E} \left[ \int_t^T \int_U v^2 \mu_s(dv) ds\right] ,
\end{equation}
where $\gamma >0$, $a=b=\sigma=\gamma=1$, $T=0.1$, $x_0=0$,
and
 $g(t,\xi ,x,y) = \gamma (x-\xi)^4/2$,
 $h(t,s,\xi ,x,y,u) = u^2/2$,
 in the framework of \eqref{cost_functional_discrete_extended}.
To the best of our knowledge, this problem admits
no analytic solution, hence we construct a numerical approximation of
its equilibrium control
based on Theorems~\ref{main theorem} and \ref{main theorem approx}
and the numerical solution $\bar{u}^{*{n},{m}}$ of Problem~\eqref{problem NM-person}.
 In Figure~\ref{fig:quartic_value_func(a)},
  we plot the value functions $J^{n,m}(0, x_0, \bar{u}^{*n,m})$
 for ${n} = 5, 10, 15, 20$ and ${m} \in \{ 1, \ldots , 20 \}$.

\begin{figure}[H]
  \centering
  \includegraphics[width=\textwidth]{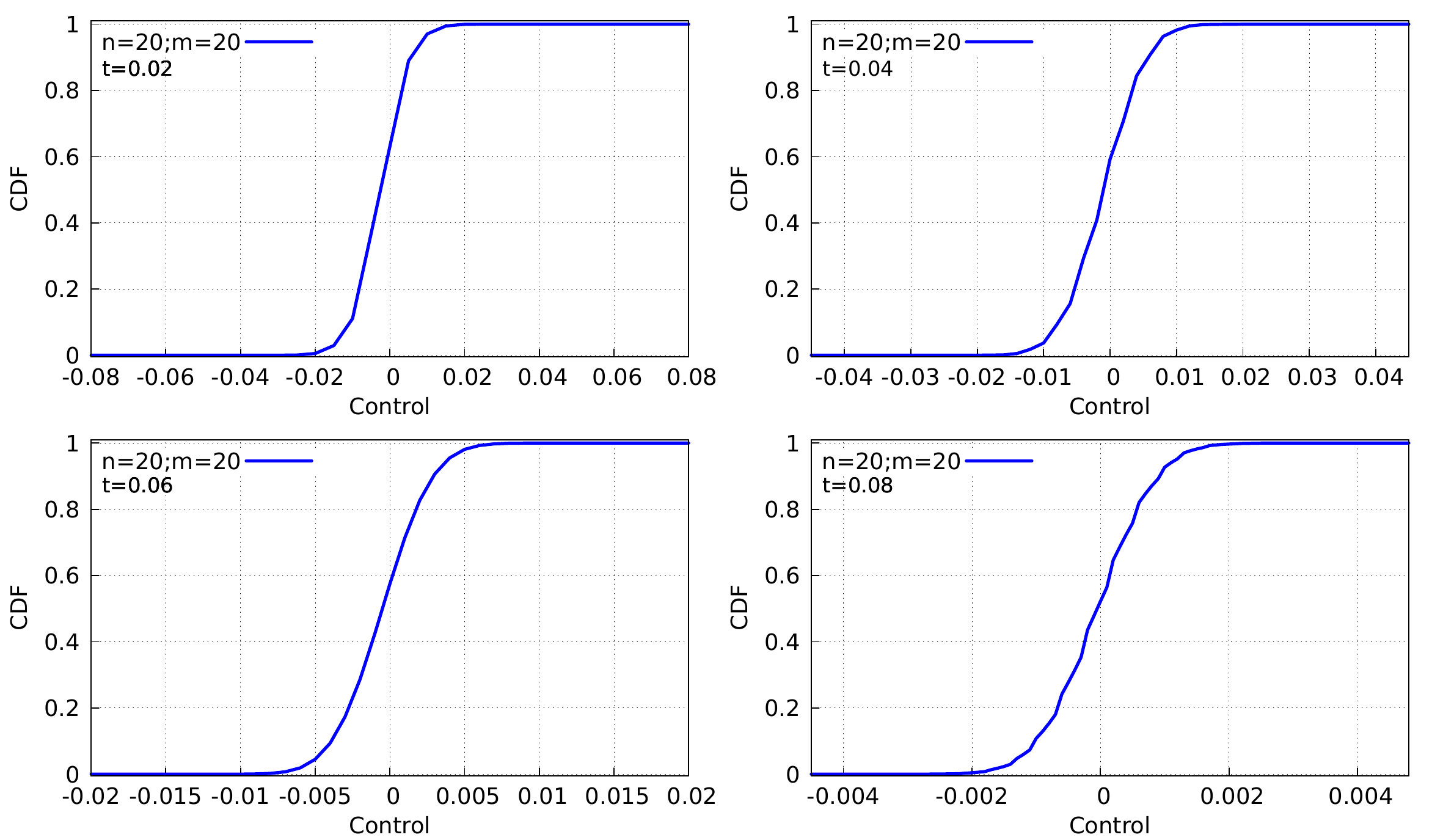} 
  \caption{CDFs at different times.}
\label{fig:quartic_cdf_diff_time}
\end{figure}

\noindent
In Figure~\ref{fig:quartic_cdf_diff_time}, we present the CDFs of
$\bar{u}^{*n,m}_t$
 with $n = m = 20$ at times $t = 0.02, 0.04, 0.06, 0.08$.

\vskip-0.2cm
\appendix
\section{Appendix} 
\noindent
The proof of Theorem~\ref{main theorem} uses Lemma~\ref{convergence implies uniform} below,
which requires the Skorokhod representation theorem
in order to construct all random variables on a single underlying probability space
as in \cite{kushner1990numerical}.
\begin{theorem}[Skorokhod representation theorem, see Theorem~6.7 in
     \cite{billingsley1999}]
  \label{skorohod theorem}
  Let $(\P_n)_{n \geq 1}$ and $\P$ be probability measures on a metric space $S$
  such that $( \P_n )_{n \geq 1}$ converges weakly to $\P$ on $S$ and
  the support of $\P$ is separable.
  Then there exist a random variable $X$ and a sequence $(X_n)_{n \geq 1}$ of random variables defined on a common probability space $(\Omega, \mathcal{F}, \P ) $, such that $\mathcal{L}(X_n) = \P_n $, $\mathcal{L}(X) = \P $, and $(X_n)_{n\geq 1}$ converges
  to $X$, $\P$-$a.s.$ on $S$.
\end{theorem}
\noindent
The following lemma was proved in
Lemma~2.4 in \cite{jacodmemin}
and Theorem~3 in \cite{valadier1994course},
and is included for completeness.
 Stable convergence of measures,
 see \cite{renyi1963},
 is defined using the test function space
 of bounded measurable functions
 $f:[0,T]\times U \to \mathbb{R}$
 such that $f(t,\cdot )$ is continuous in $u\in U$
 for all $t\in [0,T]$.
We respectively denote by ${\cal C}_b(\real )$ and
${\cal C}_b(U)$ the spaces of bounded continuous
functions on $\real$ and $U$.
\begin{lemma}\label{equivalence stable weak}
  Consider a family $( \lambda^n )_{n\geq 1} \subset \Lambda$ and $\lambda \in \Lambda$.
  The following are equivalent:
\begin{enumerate}[i)]
\item The sequence $(\lambda^n)_{n\geq 1}$ converges stably to $\lambda \in \Lambda$.
\item $\lambda^n\biggl(
    \sum\limits_{l=1}^{m} \mathbbm{1}_{A_l}(t) g_l(u)\biggr)
    \xrightarrow[n \to \infty]{} \lambda \biggl(
    \sum\limits_{l=1}^{m} \mathbbm{1}_{A_l}(t) g_l(u)\biggr)$,
    for any $m \geq 1$, any finite $\mathcal{B}([0,T])$-partition
    $\bigl\{A_1, A_2, \ldots , A_m\bigr\}$
    of $[0,T]$ and
    $g_1,\ldots , g_m \in {\cal C}_b(\real )$,
  \item The sequence $(\lambda^n)_{n\geq 1}$ converges weakly to $\lambda \in \Lambda$.
\end{enumerate}
\end{lemma}
\begin{Proof}
  As $(i) \Rightarrow (ii)$ and $(i) \Rightarrow (iii)$ are straightforward,
  we only show the following.

\noindent $(ii) \Rightarrow (i)$:
Let $f$ be a bounded measurable function $f(t,u)$
 such that $f(t,\cdot )$ is continuous in $u\in U$
 for all $t\in [0,T]$.
 By the Riesz Theorem, see \S~12.3 page~251 of \cite{royden1988real},
  the space ${\cal C}_b(U)$ is separable.
Denoting by $(c_l)_{l \geq 0}$ a countable dense subset of
${\cal C}(U)$ with respect to $\norm{\cdot }_{\infty}$,
 with $c_0\equiv 0$ and 
 letting
\[D'_{l,m}=\big\{t \in [0,T] \ : \ \norm{f(t, \cdot )-c_l}_\infty = \min\limits_{0 \leq k \leq m} \norm{f(t, \cdot )-c_k }_\infty\big\},
\quad m \geq 0, \quad l=0,1,\ldots m,
\]
 we partition $[0,T]$ into the measurable sets
 $$D_{l,m} = D'_{l,m} \setminus \bigcup\limits_{k=0}^{l-1}D'_{k,m}
 $$
 made of $t \in [0,T]$
such that $l \in \{ 0,1,\ldots , m \}$ is
 the smallest integer satisfying
 $\norm{f(t, \cdot )-c_l}_\infty = \min\limits_{0 \leq k \leq m} \norm{f(t, \cdot )-c_k}_\infty$.
 Letting $f_m(t,u) := \sum\limits_{l=0}^{m}\mathbbm{1}_{D_{l,m}}(t)c_l(u)$,
 by the denseness of $( c_l )_{l \geq 1}$ in ${\cal C}(U)$ we have
\[
\lim\limits_{m \to \infty}
\norm{f(t, \cdot )-f_m(t, \cdot )}_\infty
= 0,
\qquad
 t \in [0,T].
\]
 Since $c_0 \equiv 0$ we have
  $\min\limits_{0 \leq l \leq m} \norm{f(t, \cdot )-c_l }_\infty \leq \norm{f(t, \cdot )}_\infty$, and
 $\norm{f_m(t, \cdot )}_\infty \leq 2\norm{f(t, \cdot )}_\infty$,
 $t \in [0,T]$.
 By the uniform boundedness of $f$ and $f_m$, $m\geq 0$, 
  we have
 \begin{align}
   \nonumber 
  \sup\limits_{\lambda \in \Lambda}
  \abs{\lambda(f) - \lambda(f_m)} &= \sup\limits_{\lambda \in \Lambda}
  \left| \int_0^T \int_U f(t,v)-f_m(t,v) \lambda_t(dv) dt\right| \\
  &\leq \int_0^T \norm{f(t, \cdot )-f_m(t,\cdot )}_\infty dt
  \xrightarrow[m \to \infty]{} 0. \nonumber
\end{align}
 Therefore, for any $\varepsilon > 0$, picking $m$ such that
\[
\sup\limits_{\lambda \in \Lambda} \abs{\lambda(f) - \lambda(f_m)} < \frac{\varepsilon}{3},
\]
 and $N$ such that for all $n > N$ by $(ii)$, we have
\[
\big|\lambda^n(f_m) - \lambda(f_m)\big| < \frac{\varepsilon}{3},
\]
 hence
\begin{equation} \label{equivalence part 1}
\big|\lambda^n(f) - \lambda(f)\big| \leq \big|\lambda^n(f) - \lambda^n(f_m)\big| + \big|\lambda^n(f_m) - \lambda(f_m)\big| + \abs{\lambda(f_m) - \lambda(f)} < \varepsilon,
\end{equation}
which shows $(i)$.

\noindent
 $(iii) \Rightarrow (ii)$:
Let $f(t,u) = \sum\limits_{l=1}^{m} \mathbbm{1}_{A_l}(t) g_l(u)$ be given as in $(ii)$.
Reasoning as in \eqref{equivalence part 1}, it suffices to show that for any given $\varepsilon > 0$, we can find bounded functions $f^{(\varepsilon )} (t,u)$
 continuous in both $t\in [0,T]$ and $u \in U$, and such that
\[
\sup\limits_{\lambda \in \Lambda} \abs{\lambda(f) - \lambda(f^{(\varepsilon )} )} < \varepsilon.
\]
Denoting by $K$ the bounding constant on
$g_1,\ldots , g_m$,
by Lusin's Theorem, see e.g. Exercise~2.44 in \cite{Folland1999}, for each $\mathbbm{1}_{A_l}(t)$ we can find a closed set $F^{(\varepsilon )}_l$ such that $[0,T]\setminus F^{(\varepsilon )}_l$ has Lebesgue measure
${\rm Leb}( F^{(\varepsilon )}_l) \leq \varepsilon / (2mK)$
 and $\mathbbm{1}_{A_l}(t)$ is continuous on $F^{(\varepsilon )}_l$.
 By Tietze's extension theorem,
 see Theorem~4.16 in \cite{Folland1999},
  we can find
 a continuous extension $f^{(\varepsilon )}_l(t)$ of $\mathbbm{1}_{A_l}(t)$ from $F^{(\varepsilon )}_l$ to $[0,T]$ such that  $|f_l|$ is bounded by $1$,
 $l=1,\ldots , m$.
 Letting $f^{(\varepsilon )} (t,u) := \sum\limits_{l=1}^{m} f_l^{(\varepsilon )}(t) g_l(u)$,
 we have
\begin{align*}
  \sup\limits_{\lambda \in \Lambda} \abs{\lambda(f) - \lambda(f^{(\varepsilon )} )} &\le
  K
  \sup\limits_{\lambda \in \Lambda} \int_0^T \int_U \sum\limits_{l=1}^{m} \abs{\mathbbm{1}_{A_l}(t) - f_l(t)} \lambda_t(dv) dt \\
  &=  K \sum\limits_{l=1}^{m} \int_0^T \abs{\mathbbm{1}_{A_l}(t) - f_l(t)} dt
  \\
  &\leq 2 K \sum\limits_{l=1}^{m} {\rm Leb}( F^{(\varepsilon )}_l)
   = \varepsilon.
\end{align*}
\end{Proof}
\noindent
 The following technical lemma has been used in the proofs of
Theorem~\ref{main theorem}, Corollary~\ref{main corollary},
and Lemma~\ref{x L2 convergence proof}.
\begin{lemma}
  \label{convergence implies uniform}
  Let
  $( \mu^n )_{ n \geq 1 }\subset \mathcal{R}([0,T] )$
  be a sequence of $\Lambda$-valued relaxed controls converging
  weakly to $\mu \in \mathcal{R}( [ 0, T ] )$. Then, for any bounded random function $f:[0,T] \times U \times \Omega \to \mathbb{R}$ such that
  $f(t,\cdot ,\omega)$ is continuous
  for all $(t, \omega) \in [0,T]\times \Omega$, we have
\begin{equation}
  \lim_{n\to \infty}
  \int_{U} f(t,v,\omega) \mu^n_t(dv)
  =
  \int_{U} f(t,v,\omega) \mu_t(dv),
  \quad
  \text{a.e. } t \in [0,T], \quad \P\mbox{-}a.s.
\end{equation}
\end{lemma}
\begin{Proof}
    Since $( \mu^n )_{ n \geq 1 }$ is a sequence of random measures converging weakly to $\mu^*$, by the Skorokhod representation Theorem~\ref{skorohod theorem} there exists
  $\widetilde{\Omega}\in\mathcal{F}$ with $\P \big(\widetilde{\Omega}\big) = 1$, such that
    for all $\omega \in \widetilde{\Omega}$, $( \mu^{*n}(\omega) )_{ n \geq 1 }$ is a sequence of deterministic measures converging weakly to $\mu^*(\omega)$.
    Since the function $\mathbbm{1}_A(t)f(t,u,\omega)$ is bounded,
measurable in $t$ and continuous in $u$
for all $A \in \mathcal{B}([0,T])$,
by Lemma~\ref{equivalence stable weak} we have
\[
\lim_{n\to \infty}
\int_A \int_{U} f(t,v,\omega) \mu^n_t(\omega)(dv)dt
=
\int_A \int_{U} f(t,v,\omega) \mu_t(\omega)(dv)dt,
\qquad A \in \mathcal{B}([0,T]),
\]
 hence
\[
\lim_{n\to \infty}
\int_{U} f(t,v,\omega) \mu^n_t(\omega)(dv) =
 \int_{U} f(t,v,\omega) \mu_t(\omega)(dv), \qquad \text{a.e. } t \in [0,T], \quad \P\mbox{-}a.s.
\]
\end{Proof}
\noindent
The following lemma,
which has been used in the proofs of
Theorems~\ref{MP N person} and \ref{main theorem},
can be proved from the almost Lipschitz property of
 uniformly continuous functions.
\begin{lemma}\label{lemma: uniform and lipschitz}
  Let $X$ be a real-valued stochastic process and let
   $(X^{(\varepsilon)})_{\varepsilon \geq 0}$ be a family of real-valued stochastic processes such that for any $p \geq 1$, we have
\[
\lim\limits_{\varepsilon \downarrow 0} \sup\limits_{t \in [0,T]} \mathbb{E} \big[\big| X^{(\varepsilon)}_t - X_t\big|^{2p} ] = 0.
\]
Then, for any uniformly continuous function $f:\mathbb{R} \to \mathbb{R}$ and any $p \geq 1$, we have
\[
\lim\limits_{\varepsilon \downarrow 0} \sup\limits_{t \in [0,T]} \mathbb{E}\big[\big|f\big(X^{(\varepsilon)}_t\big) - f(X_t)\big|^{2p}] = 0.
\]
\end{lemma}
\begin{Proof}
We shall prove that for any $\varepsilon > 0$,
\[
\lim\limits_{\varepsilon \downarrow 0} \sup\limits_{t \in [0,T]} \mathbb{E}\big[\big|f\big(X^{(\varepsilon)}_t\big) - f(X_t)\big|^{2p}\big] \leq \varepsilon.
\]
 Since $f$ is uniformly continuous, for any $\rho > 0$, we can pick $K_{\rho} > 0$ such that for all $x, y \in \mathbb{R}$, we have
\[
\abs{f(x) - f(y)} \leq \rho + K_{\rho} \abs{x-y},
\]
which implies
\[
\lim\limits_{\varepsilon \downarrow 0} \sup\limits_{t \in [0,T]} \mathbb{E}\big[\big|f\big(X^{(\varepsilon)}_t\big) - f(X_t)\big|^{2p}\big] \leq 2^{2p} \big(\rho^{2p} + K_{\rho}^{2p} \lim\limits_{\varepsilon \downarrow 0} \sup\limits_{t \in [0,T]}\mathbb{E} \big[\big|X^{(\varepsilon)}_t - X_t\big|^{2p} \big]\big) = |2\rho|^{2p}.
\]
We conclude by taking $\rho = \varepsilon^{1/(2p)}/2$.
\end{Proof}

\footnotesize

\def\cprime{$'$} \def\polhk#1{\setbox0=\hbox{#1}{\ooalign{\hidewidth
  \lower1.5ex\hbox{`}\hidewidth\crcr\unhbox0}}}
  \def\polhk#1{\setbox0=\hbox{#1}{\ooalign{\hidewidth
  \lower1.5ex\hbox{`}\hidewidth\crcr\unhbox0}}} \def\cprime{$'$}

\end{document}